\documentclass[12pt]{article}

\usepackage{epsfig,graphics} 
 
\usepackage{amssymb}

\newcommand{\bo}{{\hfill\loota}} 
\newcommand{\loota}{\hbox{\enspace{\vrule height 7pt depth 0pt width 
      7pt}}}

\newcommand{\bPi}{\mbox{\boldmath $\Pi$}} 
\newcommand{\bPsi}{\mbox{\boldmath $\Psi$}} 
\newcommand{\bP}{\mbox{\boldmath $P$}}

\newcommand{\beqn}{\begin{eqnarray}} 
\newcommand{\eeqn}{\end{eqnarray}} 
\newcommand{\be}{\begin{equation}} 
\newcommand{\ee}{\end{equation}} 
\newcommand{\ba}{\begin{array}} 
\newcommand{\ea}{\end{array}} 
\newcommand{\R}{{\rm\bf R}} 
\newcommand{\C}{{\rm\bf C}} 
\newcommand{\pa}{\partial} 
 
\newcommand{\re}{\ref} 
\newcommand{\ci}{\cite} 
\newcommand{\la}{\label} 
\newcommand{\bfr}{\begin{flushright}} 
\newcommand{\efr}{\end{flushright}} 
\newcommand{\bfl}{\begin{flushleft}} 
\newcommand{\efl}{\end{flushleft}} 
\newcommand{\fr}{\frac} 
\newcommand{\ov}{\overline} 
\newcommand{\ti}{\tilde} 
\newcommand{\st}{\stackrel}

\newcommand{\si}{\sigma} \newcommand{\Si}{\Sigma} 
\newcommand{\al}{\alpha}

\newcommand{\ds}{\displaystyle}

\newcommand{\cE}{{\cal E}} \newcommand{\cF}{{\cal F}} 
\newcommand{\cH}{{\cal H}} 
 
\newcommand{\cK}{{\cal K}} 
\newcommand{\cL}{{\cal L}} 
\newcommand{\cO}{{\cal O}} 
\newcommand{\cS}{{\cal S}}\newcommand{\bS}{{\bf S}} 
 
\newcommand{\bI}{{\bf I}} 
 \newcommand{\cT}{{\cal T}} 
\newcommand{\cZ}{{\cal Z}}

\newcommand{\5}{{\hspace{0.5mm}}} 
\newcommand{\ve}{\varepsilon} 
 
\newcommand{\we}{\wedge} 
\newcommand{\de}{\delta}\newcommand{\De}{\Delta} 
 
 \newcommand{\ga}{\gamma} 
\newcommand{\om}{\omega} 
\newcommand{\Om}{\Omega} 
 
\newcommand{\na}{\nabla}

\newcommand{\lam}{\lambda} \newcommand{\ka}{\kappa} 
\newcommand{\Lam}{\Lambda} 
\newcommand{\co}{{\rm const}} 
\newcommand{\supp}{{\rm supp\5}}

\newcommand{\brr}{{|\kern-.15em|\kern-.15em|\kern-.15em}\,} 
\newcommand{\ddd}{\st{.\kern-.07em.\kern-.07em.}} 
\def\N{{\rm I\kern-.1567em N}}                              
\def\R{{\rm I\kern-.1567em R}}                              
\def\C{{\rm C\kern-4.7pt                                    
\vrule height 7.7pt width 0.4pt depth -0.5pt \phantom {.}}} 
\def\Z  {{\sf Z\kern-4.5pt Z}}                                
\def\Re {{\rm Re\, }}                                       
\def\Im {{\rm Im\,}}                                        
 
\begin{document} 
 
\renewcommand{\theequation}{\thesection.\arabic{equation}} 
\newtheorem{theorem}{Theorem}[section] 
\renewcommand{\thetheorem}{\arabic{section}.\arabic{theorem}} 
\newtheorem{definition}[theorem]{Definition} 
\newtheorem{deflem}[theorem]{Definition and Lemma} 
\newtheorem{lemma}[theorem]{Lemma} 
\newtheorem{example}[theorem]{Example} 
\newtheorem{remark}[theorem]{Remark} 
\newtheorem{remarks}[theorem]{Remarks} 
\newtheorem{cor}[theorem]{Corollary} 
\newtheorem{pro}[theorem]{Proposition}

\newcommand{\bd}{\begin{definition}} 
 \newcommand{\ed}{\end{definition}} 
\newcommand{\bt}{\begin{theorem}} 
 \newcommand{\et}{\end{theorem}} 
\newcommand{\bqt}{\begin{qtheorem}} 
 \newcommand{\eqt}{\end{qtheorem}}

\newcommand{\bp}{\begin{pro}} 
 \newcommand{\ep}{\end{pro}} 
 
\newcommand{\bl}{\begin{lemma}} 
 \newcommand{\el}{\end{lemma}} 
\newcommand{\bc}{\begin{cor}} 
 \newcommand{\ec}{\end{cor}} 
 
\newcommand{\bex}{\begin{example}} 
 \newcommand{\eex}{\end{example}} 
\newcommand{\bexs}{\begin{examples}} 
 \newcommand{\eexs}{\end{examples}}

\newcommand{\bexe}{\begin{exercice}} 
 \newcommand{\eexe}{\end{exercice}}

\newcommand{\br}{\begin{remark} } 
 \newcommand{\er}{\end{remark}} 
\newcommand{\brs}{\begin{remarks}} 
 \newcommand{\ers}{\end{remarks}}

\newcommand{\pru}{{\bf Proof~~}}

\begin{titlepage} 
 
\begin{center} 
{\Large\bf 
On Scattering of Solitons for the Klein-Gordon\\ 
~\\ 
Equation Coupled to a Particle 
}\\ 
\vspace{2cm} 
{\large Valery Imaikin} 
\footnote{Supported partly by Austrian Science Foundation 
(FWF) Project P19138-N13, by research grants of DFG 
(436 RUS 113/615/0-1(R)) and RFBR (01-01-04002).} 
\medskip\\ 
{\it Wolfgang Pauli Institute\\ 
 c/o Faculty of Mathematics of 
Vienna University\\ 
 Nordbergstrasse 15, 1090 Vienna, Austria}\\ 
email: valery.imaikin@univie.ac.at\bigskip\\ 
{\large Alexander Komech}$^1\,$ 
\footnote{ 
On leave Department Mechanics and Mathematics 
of Moscow State University. 
Supported partly  by Austrian Science Foundation 
(FWF) Project P19138-N13,
by 
Max-Planck Institute of Mathematics in the Sciences (Leipzig), 
 and Wolfgang Pauli Institute of Vienna University.} 
\medskip\\ 
{\it Faculty of Mathematics of 
Vienna University\\ 
Nordbergstrasse 15, 
1090 Vienna, Austria 
}\\ 
email: alexander.komech@univie.ac.at\bigskip\\ 
{\large Boris Vainberg}\footnote{Supported partially by the NSF grant 
DMS-0405927}\medskip\\ 
{\it Department of Mathematics and Statistics\\ UNC at Charlotte, 
Charlotte, NC 28223, USA}\\ 
email: brvainbe@mail.uncc.edu\bigskip\\ 
 
\end{center} 
 
\vspace{2cm} 
 
\begin{abstract} 
 
We establish the long time soliton asymptotics for the translation invariant
nonlinear system consisting 
of the Klein-Gordon equation coupled to a charged relativistic
particle. The coupled system has a six dimensional invariant manifold of the soliton
solutions. We show that in the large time approximation any finite energy solution,
with the initial state close to the solitary manifold, is a sum of a
soliton and a dispersive wave which is a solution of the free Klein-Gordon
equation. It is assumed that the charge density satisfies the Wiener
condition which is a version of the ``Fermi Golden Rule''. The proof is based
on an extension of the general strategy 
introduced
by 
Soffer and Weinstein,
Buslaev and Perelman,  and others:
symplectic projection in Hilbert space onto the solitary manifold, modulation 
equations for the parameters of the projection, and decay of the transversal component.

\end{abstract} 
\end{titlepage} 
 
 
\setcounter{equation}{0} 
 
\section{Introduction} 
 
Our paper concerns the problem of nonlinear 
field-particle interaction. A charged particle radiates a field 
which acts back on the particle. This interaction is  
responsible for some crucial features of the dynamics: asymptotically 
 uniform motion and stability against small perturbations 
of the particle, increase of the particle's mass and others 
(see \ci{Abr,Dir, Lor,Sp}). 
The problem has many different appearances: a classical particle 
coupled to a scalar or Maxwell field, and coupled Maxwell-Schr\"odinger 
 or Maxwell-Dirac equations, general translation invariant   
nonlinear hyperbolic PDEs. 
In all cases the goal is to reveal  
the distinguished role of the soliton solutions,  
i.e. traveling wave solutions of finite energy. 
Let us note that the existence of the soliton solutions is proved for  
nonlinear Klein-Gordon equations with a general nonlinear term 
\ci{BL}, and for the coupled  Maxwell-Dirac equations \ci{EGS}. 
 
One of the main goals of  
a mathematical investigation 
is to study soliton type asymptotics and 
asymptotic stability of soliton solutions to the equations.  
First results in this direction have  
been discovered for the KdV equation and 
other {\it completely integrable equations}.  
For the KdV equation, any 
solution with sufficiently smooth and rapidly decaying initial data 
converges to a finite sum of soliton solutions moving to the right, and a 
dispersive wave moving to the left. A complete survey and proofs can 
be found in \ci{EH}.

For nonintegrable equations, the long time  
convergence of the solution 
to a soliton part and dispersive wave was obtained first  
by Soffer and Weinstein in the context of   
 $U(1)$-invariant  Schr\"odinger equation \ci{SW88,SW90,SW92}. 
The extension to translation invariant equations was obtained by  
 Buslaev and Perelman \ci{BP1,BP2} for the 1D Schr\"odinger equation, 
and by Miller, Pego and  
Weinstein for the 1D modified KdV and RLW equations, \ci{MW96,PW92,PW94}. 
The techniques introduced by Weinstein \ci{W85} 
play a fundamental role 
in the proofs of all these results.

In \ci{BP1,BP2} the long time convergence is obtained for the 1D 
translation invariant 
and $U(1)$-invariant  
nonlinear Schr\"odinger equation. It is shown there that 
the following asymptotics hold 
for any 
finite-energy solution $\psi(x,t)$ with initial data close to a 
 soliton $\psi_{v_0}(x-v_0t-a_0)e^{i\om_0t}$: 
\be\la{Bus} 
\psi(x,t)=\psi_{v_\pm}(x-v_\pm t-a_\pm)e^{i\om_\pm 
t}+W_0(t)\psi_{\pm}+ r_\pm(x,t),\,\,\,t\to\pm\infty. 
\ee  
Here the 
first term on the right hand side is a soliton with parameters 
$v_\pm$, $a_\pm$, $\om_\pm$ close to $v_0$, $a_0$, $\om_0$, the 
function  $W_0(t)\psi_{\pm}$ is a dispersive wave which is a 
solution to the free Schr\"odinger equation, 
 and the remainder $r_\pm(x,t)$ converges to zero in the global $L^2$-norm. 
Recently Cuccagna extended the asymptotics 
(\re{Bus}) 
 to nD Schr\"odinger equations with 
$n\ge 3$, \ci{Cu, Cu03}. 
 
We establish the asymptotics similar to 
(\re{Bus}) for 
a scalar real-valued Klein-Gordon field $\psi (x)$ in $\R^3$ 
coupled to a relativistic particle with position $q$ and momentum $p$ 
governed by 
\beqn\la{system0} 
\ba{lll} 
\dot{\psi}(x,t)=\pi (x,t), & 
\dot{\pi}(x,t)=\Delta \psi (x,t) 
-m^2\psi(x,t)-\rho (x-q(t)), & x\in\R^3, 
\\ 
~ &  \\ 
\dot{q}(t)=p(t)/\sqrt{1+p^2(t)}, & \dot{p}(t)=\displaystyle\int 
\psi 
(x,t)\,\nabla \rho (x-q(t))dx, & 
\ea 
\eeqn 
 where $m>0$ (the case $m=0$ is degenerate and will be considered 
elsewhere). 
This is a Hamiltonian system with the Hamiltonian functional 
\be\la{hamilq0} 
{\cal H}(\psi ,\pi ,q,p)=\frac 12\int\Big(|\pi(x)|^2+ 
|\nabla \psi (x)|^2+m^2|\psi(x)|^2\Big)dx+\int \psi (x)\rho (x-q)dx+ 
\sqrt{1+p^2}. 
\ee 
The first two equations for the fields are equivalent to the 
 Klein-Gordon equation with the source $\rho(x-q)$. 
The form of the last two equations in 
(\re{system0}) 
is determined by the choice of the relativistic 
kinetic energy $ \sqrt{1+p^2}$ in (\re{hamilq0}). 
Nevertheless, the system (\re{system0}) is not relativistic invariant. 
 
We have set the maximal speed of the particle equal 
to one, which is the speed of wave propagation. 
This is in agreement with the principles of special relativity. 
Let us also note that the first two equations 
of  (\re{system0}) admit the soliton solutions of finite energy, 
$\psi_v(x-vt-a),\pi_v(x-vt-a)$, if and only if $|v|<1$.

The case of a point particle corresponds to $\rho (x)=\de(x)$ and then 
the interaction term in the Hamiltonian is simply 
$\psi (q)$. 
However, in this case the Hamiltonian is unbounded from below which 
leads to the ill-posedness of the problem, 
also known as ultraviolet divergence. 
Therefore we smooth the coupling by the function $\rho (x)$ 
following the 
``extended electron'' 
strategy proposed by M.~Abraham \ci{Abr} for  
charges coupled to 
the Maxwell field. 
 In analogy to the Maxwell-Lorentz equations we call $\rho$ the 
``charge distribution''. 
Let us write the system (\re{system0}) 
as 
\be\la{HPDE} 
\dot Y(t)=F(Y(t)),~~~~~t\in\R, 
\ee 
where $Y(t):=(\psi(x,t),\pi(x,t),q(t),p(t))$ 
(below we always deal with column vectors but often write them as 
row vectors).  
The system (\re{system0}) is translation-invariant 
and admits the soliton solutions 
\be\la{sosol} 
Y_{a,v}(t)=(\psi_v(x-vt-a),\pi_v(x-vt-a), vt+a,p_v), 
~~~~~~~p_v=v/{\sqrt{1-v^2}} 
\ee 
for all $a,v\in\R^3$ with $|v|<1$ (see (\re{stfch}), (\re{sol})), 
where the functions $\psi_v$, $\pi_v$ decay exponentially for $m>0$ 
(main difficulty of the case $m=0$ is provided by very slow decay of the functions). 
The states $S_{a,v}:=Y_{a,v}(0)$ form the solitary manifold 
\be\la{soman} 
{\cal S}:=\{ S_{a,v}: a,v\in\R^3, |v|<1 \}. 
\ee 
Our main result is the soliton-type asymptotics 
of type (\re{Bus}) for $t\to\pm\infty$, 
\be\la{Si} 
(\psi(x,t), \pi(x,t))\sim 
(\psi_{v_\pm}(x-v_\pm t-a_\pm), \pi_{v_\pm}(x-v_\pm t-a_\pm)) 
+W_0(t)\bPsi_\pm 
\ee 
for solutions 
to  (\re{system0}) 
with initial data close 
to the solitary manifold ${\cal S}$. 
Here $W_0(t)$ is the dynamical group of the 
free Klein-Gordon equation, 
$\bPsi_\pm$ are the corresponding asymptotic scattering states, 
and the remainder converges to zero 
in the global energy norm, i.e. in the norm of the 
Sobolev space 
$H^1(\R^3)\oplus L^2(\R^3)$. 
For the particle 
trajectory we prove that 
\be\la{qqi} 
\dot q(t)\to v_\pm, 
~~~~q(t) \sim v_\pm t+a_\pm. 
\ee 
The results are established under the following conditions 
on the charge distribution: $\rho$ is a real valued function of the 
Sobolev class $H^2(\R^3)$, compactly supported, 
and spherically symmetric, i.e. 
\be\la{ro} 
\rho,\na\rho,\na\na\rho\in L^2({\R}^3),~~~~~\,\,\,\,\,\,\quad\rho (x)=0 
\,\,\,\mbox{for}\,\,\,|x|\ge R_\rho, \,~~~~~~~~~~\rho (x)=\rho_1(|x|). 
\ee 
We require that all ``modes'' of the wave field are 
coupled to the particle, which is formalized by the Wiener 
condition 
\be\la{W} 
\hat\rho(k)=(2\pi)^{-3/2}\int\limits \, e^{ik x}\rho(x)dx\not=0 
\mbox{ \,\,\,for\,\,all\,\, }k\in\R^3\,. 
\ee 
It is an analogue of the ``Fermi Golden Rule'' 
\ci{BP3,BS,Cu,Cu03,Sig,SW3,SW4}: 
the coupling term  $\rho(x-q)$ is not orthogonal 
to the eigenfunctions $e^{ikx}$ 
of the continuous spectrum of the linear part 
of the   equation. 
As we will see, the Wiener condition (\re{W}) is very essential 
for our asymptotic analysis (see Remark \re{reW}). 
Generic examples of the coupling function $\rho$ satisfying (\re{ro}) and 
 (\re{W}) are given in \ci{KSK}. 
\br 
{\rm 
Physically, the Wiener condition means the strong coupling 
of the particle to the field  
which leads to 
{\it radiation} 
of the particle.  
This radiation results 
in the 
relaxation of the acceleration $\ddot q(t)\to 0$, $t\to\pm\infty$ 
which provides 
the asymptotics  (\re{Si}) and (\re{qqi}). 
Note that the soliton solutions do not radiate, and 
the radiation  
of the particle 
manifests itself in the  
decay of the 
deviation of the solution 
from the solitary manifold (see (\re{decnonlin}) below). 
} 
\er

The problem under investigation was studied earlier in  
the following 
two different 
situations A and B: 
 
\noindent A. The asymptotics 
\be\la{as0} 
\dot q(t)\to v_\pm,~~~~~~ 
(\psi(x,t),\pi(x,t))\sim(\psi_{v_\pm}(x-q(t)),\pi_{v_\pm}(x-q(t))) 
\ee 
were proved in \ci{KS} in the case $m=0$, under the Wiener condition (\re{W}), 
for all finite energy solutions, without the assumption that the initial 
data are close to ${\cal S}$. This means that the solitary manifold 
is a {\it global attractor} for the equations (\re{system0}). However, 
the asymptotics (\re{as0}) were established only in  {\it local energy 
semi-norms} centered at the particle position $q(t)$. This means that 
the remainder in (\re{as0}) may contain a dispersive term, 
similar to 
the middle term in the right hand side of (\re{Si}), whose energy 
 radiates to infinity as $t\to\pm\infty$ but does not converge to zero. 
  A similar result is established 
in  \ci{IKMml} for coupled Maxwell-Lorentz equations. 
\medskip 
 
\noindent B. The asymptotics (\re{as0}),  
and an analogue of the asymptotics    
(\re{Si}) in the 
{\it global 
energy norm}, were established in \ci{IKSs} (resp., \ci{IKM}) 
also for all finite energy solutions, 
in the case $m=0$ (resp., $m>0$), 
under the smallness 
condition on the coupling function, 
 $\Vert\rho\Vert_{L^2(\R^3)}\ll1$.  
The similar results are established 
in  \ci{IKSm,Sp} (resp., \ci{IKSr}) for the coupled Maxwell-Lorentz equations 
with a moving (resp., rotating) charge. 
Let us stress that the asymptotics   (\re{qqi}) 
for the position was missing in the previous work.

\medskip 
 
\noindent 
 
Let us comment on main difficulties in proving the 
 asymptotic stability of the invariant manifold $S$ and justifying 
(\re{Si}), (\re{qqi}). 
The method of \ci{IKMml,KS} is based on the Wiener Tauberian 
Theorem, 
hence cannot provide a rate of convergence  
in the velocity asymptotics 
of (\re{qqi}) which is needed to prove (\re{Si}) and the position asymptotics  
of (\re{qqi}). 
Also the methods of \ci{IKM,IKSm,IKSs,IKSr} are applicable only 
for a small coupling function $\rho(x)$, 
and do not provide  the position asymptotics  in (\re{qqi}).

Our approach develops a general  
strategy introduced in  \ci{BP1,BP2,PW92,PW94} 
for proving the asymptotic 
stability of the invariant solitary manifold $\cS$. 
The strategy originates from the 
techniques in \ci{W85}  and their  developments in 
\ci{SW88, SW90,SW92} in the context of the $U(1)$-invariant Schr\"odinger equation. 
The 
approach uses the symplectic geometry methods for the Hamiltonian 
systems in Hilbert spaces and spectral theory of nonselfadjoint 
operators. 
 
The invariant manifolds arise automatically for 
equations 
with a symmetry Lie group \ci{BL,EGS,GSS}. 
In particular, our system (\re{system0}) is invariant 
under translations in $\R^3$. 
The asymptotic stability  
of the solitary manifold 
is studied by a linearization of 
the dynamics (\re{HPDE}). 
The linearization will be made along a special curve on 
the solitary manifold,  $S(t)$,  
which is the symplectic orthogonal projection 
of the solution.  
Then the linearized equation reads 
\be\la{lind} 
\dot X(t)=A(t)X(t),~~~~~~t\in\R, 
\ee 
where the operator $A(t)$ 
corresponds to the linearization at the soliton $S(t)$. 
Furthermore, we 
consider the ``frozen'' linearized equation (\re{lind}) with 
$A(t_1)$ instead of $A(t)$. The operator $A(t_1)$ 
has zero eigenvalue, and the frozen linearized equation 
admits secular solutions  linear in $t$ (see (\re{secs})). 
The existence of these 
runaway 
solutions 
prohibits the direct 
application of the Liapunov strategy and 
is responsible for the 
instability of the nonlinear dynamics along the manifold 
$\cS$. 
 One crucial observation is that 
the linearized equation is stable in 
the 
{\it symplectic orthogonal complement} 
to the tangent space $\cT_S$. 
The complement is invariant under the linearized dynamics 
since 
the linearized  dynamics is Hamiltonian and 
leaves the symplectic structure 
invariant. 
\bigskip\bigskip 
 
Our proofs are based on a suitable extension of the methods in 
\ci{BP1,BP2,PW92,PW94}. 
Let us comment on the main steps. 
\medskip 
 
I. First, 
we construct the symplectic orthogonal projection $S(t)=\bPi Y(t)$ 
of the trajectory 
$Y(t)$ onto the solitary manifold $\cS$. This means that 
$S(t)\in\cS$, and 
the complement vector $Z(t):=Y(t)-S(t)$ is symplectic orthogonal to 
the tangent space $\cT_{S(t)}$ for every $t\in\R$: 
\be\la{sio} 
Z(t)\nmid \cT_{S(t)},~~~~~~~~~~~~~~~~t\in\R. 
\ee 
So, we get the splitting 
$Y(t)=S(t)+Z(t)$ and we linearize the dynamics in the 
{\it transversal component} $Z(t)$ along the trajectory.

 The  {\it soliton component} $S(t)=S_{b(t), v(t)}$ 
satisfies 
a {\it modulation equation}. 
Namely, in the parametrization $\xi(t)=(c(t),v(t))$ 
with $c(t):=b(t)-\ds\int_0^t v(s)ds$, we have 
\be\la{modp} 
\dot \xi(t)=N_1(\xi(t), Z(t)),~~~~~~~~~~~ 
|N_1(\xi(t), Z(t))|\le C\Vert Z(t)\Vert_{-\beta}^2\,, 
\ee 
where 
$\Vert\cdot\Vert_{-\beta}$ 
stands for  an appropriate 
weighted Sobolev norm. 
 
On the other hand, 
the transversal component  satisfies 
the {\it transversal equation} 
\be\la{nlna} 
\dot Z(t)=A(t)Z(t)+N_2(S(t),Z(t)), 
\ee 
where $A(t)=A_{S(t)}$, and $N_2(S(t),Z(t))$ is a nonlinear part: 
\be\la{Np} 
\Vert N_2(S(t),Z(t))\Vert_{\beta}\le C\Vert Z(t)\Vert_{-\beta}^2, 
\ee 
where $\Vert\cdot\Vert_{\beta}$ is defined 
similarly to $\Vert\cdot\Vert_{-\beta}$. 
Let us note that the bound (\re{Np}) is not a direct consequence 
of the linearization, since the function $S(t)$ generally 
is not a solution of (\re{HPDE}). The modulation equation and the bound 
(\re{modp}) 
play a crucial role in the proof of (\re{Np}).

\medskip 
 
II. The linearized dynamics 
(\re{lind}) is nonautonomous. First, 
let us fix $t=t_1$ in $A(t)$ and consider the corresponding 
``frozen'' 
linear autonomous equation with $A(t_1)$ instead of $A(t)$. 
We prove the decay 
\be\la{declin} 
\Vert X(t)\Vert_{-\beta}\le \ds\fr{C\Vert X(0)\Vert_{\beta}}{(1+|t|)^{3/2}}, 
~~~~~~~t\in\R 
\ee 
of the solutions $X(t)$ to the frozen equation 
for any $X(0)\in \cZ_{S_1}$ where $S_1:=S(t_1)$, 
and $\cZ_{S_1}$ 
is the space of vectors $X$ which are symplectic orthogonal to 
the tangent space $\cT_{S_1}$. 
Let us stress that the decay holds only 
for the solutions 
symplectic orthogonal to the tangent space. 
Basically, the reason of the decay is the fact that the spectrum 
of the generator $A(t_1)$ restricted to the space $\cZ_{S_1}$ is 
purely continuous. 
\medskip 
 
III. We combine the decay (\re{declin}) 
with the bound (\re{modp}) 
through the 
nonlinear equation (\re{nlna}). 
This 
gives the time decay of the transversal 
component 
\be\la{decnonlin} 
\Vert Z(t)\Vert_{-\beta}\le 
\ds\fr{C(\Vert Z(0)\Vert_{\beta})}{(1+|t|)^{3/2}},~~~~~~t\in\R, 
\ee 
if the norm $\Vert Z(0)\Vert_{\beta}$ is sufficiently small. 
 One of the main difficulties in proving the decay 
(\re{decnonlin}) 
 is the non-autonomous 
character of the linear part of 
 (\re{nlna}). 
 We deduce the decay 
from the equation 
(\re{nlna}) written in the ``frozen'' form 
\be\la{nlnaf} 
\dot Z(t)=A(t_1)Z(t)+[A(t)-A(t_1)]Z(t)+ 
N_2(S(t),Z(t)),~~~~~~~~~~~~0\le t<t_1, 
\ee 
with arbitrary large $t_1>0$.

IV. The decay (\re{decnonlin}) implies the soliton asymptotics 
(\re{Si}) and (\re{qqi}) by the known techniques of 
scattering theory. 
 
\brs 
{\rm 
i) 
The asymptotic stability 
of the solitary manifold $\cS$ is 
caused by 
the {\it radiation of energy to infinity} 
which appears as the {\it local energy decay} 
for the transversal component, 
(\re{decnonlin}). 
\smallskip\\ 
ii) 
The asymptotics  (\re{Si}) can be interpreted as the collision 
of the incident soliton,  with a trajectory $v_-t+a_-$, 
with an incident wave $W_0(t)\bPsi_-$, which results in an 
outgoing soliton  with a new trajectory $v_+t+a_+$, 
and a new outgoing wave $W_0(t)\bPsi_+$. The collision process 
can be represented 
by the diagram of Fig. 1. 
It 
suggests to introduce the (nonlinear) {\it scattering 
operator} 
\be\la{sco} 
\bS:~(v_-,a_-,\bPsi_-)\mapsto(v_+,a_+,\bPsi_+). 
\ee 
However, the domain of the operator 
is an open problem as well as the question on its 
{\it asymptotic completeness} 
(i.e. on its range). 
} 
\ers 
\begin{figure}[htbp] 
\begin{center} 
\includegraphics[width=0.8\columnwidth]{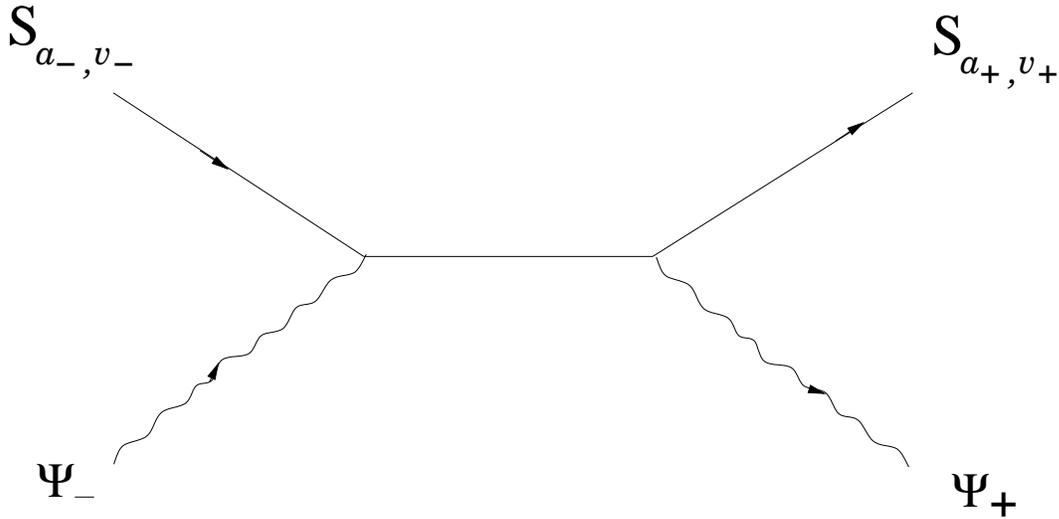} 
\caption{Wave -- particle scattering.} 
\label{fig-1} 
\end{center} 
\end{figure} 
 
\brs 
{\rm 
i) The strategy of \ci{BP1,BP2,PW92,PW94} was further developed 
in the papers \ci{BP3,BS,Cu,Cu03,MW96,SW3,SW4, SW05}. 
Let us stress that 
these papers 
contain several assumptions on the discrete and continuous 
spectrum 
of the linearized problem. In our case a complete 
investigation of the spectrum of the linearized problem 
is given under the Wiener condition 
and there is no need for any a priori spectral assumptions. 
\smallskip\\ 
ii) 
Note that the Wiener condition is indispensable for our proof 
of the decay (\re{declin}),  
but only  
in the proof of Lemma \re{lW}. 
Otherwise we use 
only the fact that the coupling function $\rho(x)$ 
is not identically zero. 
The other assumptions on $\rho$ 
can be weakened: 
the spherical symmetry 
is not necessary, and 
one can assume also that $\rho$ 
belongs to a weighted Sobolev space rather than having a compact support. 
} 
\ers

\noindent 
 
Our paper is organized as follows. In Section 2, we formulate the main 
result. In  Section 3, we introduce the symplectic projection onto the 
solitary manifold. The linearized equation is defined in section 4. In 
Section 6, we split the dynamics in two components: 
along the solitary manifold and in the transversal directions, and we 
justify the estimate (\re{modp}) 
concerning the tangential component. The time 
decay of the transversal component is established in sections 7 - 10 
under an assumption on the time decay of the linearized dynamics. In 
Section 11, we prove the main result. Sections 12 - 18 fill the 
gap concerning the time decay of the linearized dynamics. 
In Appendices A and B we collect some routine calculations. 
\bigskip\\ 
{\bf Acknowledgments} 
The authors thank V.Buslaev 
for numerous lectures on his results and 
fruitful discussions, and E.Kopylova for 
reading the paper and making many useful remarks.


\setcounter{equation}{0} 
 
\section{Main Results} 
 
 
\subsection{Existence of Dynamics} 
 
 To formulate our results precisely, we need some definitions. We introduce 
a suitable phase space for the Cauchy problem corresponding to 
(\ref{system0}) and (\ref{hamilq0}). 
Let $H^0=L^2$ be the real Hilbert space $L^2({\R}^3)$ with scalar product 
$\langle\cdot,\cdot\rangle$ and 
norm $\Vert\cdot\Vert_{L^2}$, and let $H^1$ be the 
 Sobolev space 
$H^1=\{\psi\in L^2:\,|\nabla\psi|\in L^2\}$ 
 with the norm 
$\Vert\psi\Vert_{H^1}=\Vert\nabla\psi\Vert_{L^2}+\Vert\psi\Vert_{L^2}$. 
Let us introduce also the weighted Sobolev spaces $H^s_{\alpha}$, $s=0,1$, 
$\alpha\in\R$ with the norms $\Vert\psi\Vert_{s,\alpha}:= 
\Vert(1+|x|)^{\alpha}\psi\Vert_{H^s}$. 
 
\begin{definition} 
 i) The phase space ${\cal E}$ is the real Hilbert space 
$H^1\oplus L^2\oplus {\R}^3\oplus {\R}^3$ of states 
$Y=(\psi ,\pi ,q,p)$ with the finite norm 
$$ 
\Vert Y\Vert_{\cal E}=\Vert \psi \Vert_{H^1} + 
\Vert\pi \Vert_{L^2}+|q|+|p|. 
$$ 
ii) ${\cal E}_{\alpha}$ is the space 
$H^1_{\alpha}\oplus H^0_{\alpha}\oplus {\R}^3\oplus {\R}^3$ 
with the norm 
\be\la{alfa} 
\Vert Y\Vert_{\alpha}=\Vert \,Y\Vert_{{\cal E}_{\alpha}}= 
\Vert \psi \Vert_{1,\alpha} +\Vert\pi \Vert_{0,\alpha}+|q|+|p|. 
\ee 
 iii) ${\cal F}$ is the space $H^1\oplus L^2$ of fields 
$F =(\psi ,\pi )$ with the finite norm 
$$ 
\Vert F \Vert_{\cal F}=\Vert \psi \Vert_{H^1} +\Vert\pi \Vert_{L^2}. 
$$ 
Similarly, ${\cal F}_{\alpha}$ is the space 
$H^1_{\alpha}\oplus H^0_{\alpha}$ 
with the norm 
\be\la{Falfa} 
\Vert \,F\Vert_{\alpha}=\Vert \,F\Vert_{{\cal F}_{\alpha}}= 
\Vert \psi \Vert_{1,\alpha} +\Vert\pi \Vert_{0,\alpha}. 
\ee 
\end{definition} 
Note that we use the same notation for the norms in the space 
 $\mathcal{F}_{\alpha}$ as in the space 
$\mathcal{E}_{\alpha}$ defined in (\ref{alfa}). We hope it will not 
create misunderstandings since  $\mathcal{F}_{\alpha}$ is 
equivalent to the subspace of $\mathcal{E}_{\alpha}$ which 
consists of elements of $\mathcal{E}_{\alpha}$ with zero vector 
components : $q=p=0$. It will be always clear from the context if 
we deal with fields only, and therefore with the space 
$\mathcal{F}_{\alpha}$, or with 
fields-particles, and therefore with elements of the space 
$\mathcal{E}_{\alpha}$.

 We consider the Cauchy problem for the Hamilton system (\re{system0}) 
 which we write as 
 \be\la{WP2.1} 
 \dot Y(t)=F(Y(t)),\quad t\in\R:\quad Y(0)=Y_0. 
 \ee 
 Here $Y(t)=(\psi(t), \pi(t), q(t), p(t))$, $Y_0=(\psi_0, \pi_0, q_0, p_0) 
 $, and all 
derivatives are understood in the sense of distributions.

 \begin{pro}\la{WPexistence} {\rm \ci{IKM}} Let (\re{ro}) hold. 
 Then\\ 
 (i) For every $Y_0\in {\cal E}$, 
the Cauchy problem (\re{WP2.1}) has a unique 
 solution $Y(t)\in C(\R, {\cal E})$.\\ 
 (ii) For every $t\in\R$, the map $U(t): Y_0\mapsto Y(t)$ is continuous on 
 ${\cal E}$.\\ 
 (iii) 
The  energy is conserved, i.e. 
 \be\la{2.4} 
{\cal H}(Y(t))= {\cal H}(Y_0),\,\,\,\,\,t\in\R, 
 \ee 
 and the velocity is bounded, 
 \be\la{WP2.1'} 
 |\dot q(t)|\leq \ov v<1,~~~~t\in\R, 
 \ee 
where $\ov v=\ov v(Y_0)$. 
 \end{pro} 
The proof is based on a priori estimates provided by the fact that 
the Hamilton functional  (\re{hamilq0}) is bounded from below. 
The latter follows from the bounds 
\be\la{hbound} 
-\fr{1}{2m^2}\Vert\rho\Vert^2_{L^2} 
\le\fr{m^2}{2}\Vert\psi\Vert^2_{L^2}+ 
\langle\psi,\rho(\cdot-q)\rangle\le\fr{m^2+1}{2}\Vert\psi\Vert^2_{L^2}+ 
\fr{1}{2}\Vert\rho\Vert^2_{L^2}, 
\ee 
which imply also that ${\cal E}$ is the space of finite energy states.

 
\subsection{Solitary Manifold and Main Result} 
 
Let us compute the  solitons (\re{sosol}). 
The substitution to (\re{system0}) gives the following stationary equations, 
\be\la{stfch} 
\left.\ba{rclrcl} 
-v\cdot\na\psi_v(y)&=&\pi_v(y)~,&-v \cdot\na\pi_v(y)&=&\De\psi_v(y)- 
m^2\psi_v(y)-\rho(y) 
\\\\ 
v&=&\ds\fr{p_v}{\sqrt{1+p_v^2}}~,&0&=&-\ds\int\na\psi_v(y)\rho(y)\,dy 
\ea\right| 
\ee 
Then the first two equations imply 
\be\la{Lpv} 
\Lam\psi_v(y):=[-\De+m^2+(v\cdot\na)^2]\psi_v(y)=-\rho(y),~~~~~~~~y\in\R^3. 
\ee 
For $|v|<1$ the operator $\Lam$ is an isomorphism  $H^4(\R^3)\to H^2(\R^3)$. 
Hence  (\re{ro}) implies that 
\be\la{Lpvy} 
\psi_v(y)=-\Lam^{-1}\rho(y)\in H^4(\R^3). 
\ee 
If $v$ is given and $|v|<1$, then $p_v$ can be found from the 
third equation of (\re{stfch}). Further, functions $\rho$ and 
$\psi_v$ are even 
by (\re{ro}). Thus, $\na\psi_v$ is odd and the last equation of 
(\re{stfch}) holds. Hence, the soliton 
solution 
(\re{sosol}) 
exists and defined uniquely for any couple $(a,v)$ with $|v|<1$. 
 
The function $\psi_v$ can be computed by the Fourier transform. 
The soliton 
is given by 
the formulas 
\be\la{sol} 
 \left.\ba{rcl} 
\psi_v(x) &=& \ds 
 -\fr \ga {4\pi} \int 
\fr {e^{-m|\ga(y-x)_\|+(y-x)_\bot|}\rho(y)d^3y} 
{|\ga(y-x)_\|+(y-x)_\bot|} 
\\\\ 
\pi_v(x) &=&-v\cdot \na\psi_v(x) , \,\,\,\,\,\, p_v= \ga v 
=\ds\fr v{\sqrt{1-v^2}} 
\ea\right| 
\ee 
Here we set 
$\ga=1/\sqrt{1-v^2}$ and 
$x=x_\Vert+x_\bot$, where 
$x_\Vert\Vert v$ 
and $x_\bot\bot v$ for $x\in\R^3$. 
 
Let us denote by  $V:=\{v\in\R^3:|v|<1\}$. 
\begin{definition} 
A soliton state is $S(\si):=(\psi_v(x-b),\pi_v(x-b),b,p_v)$, where 
 $\si:=(b,v)$ with 
 $b\in\R^3$ and $v\in V$. 
\end{definition} 
Obviously, the soliton solution admits the representation $S(\si(t))$, where 
\be\la{sigma} 
\si(t)=(b(t),v(t))=(vt+a,v). 
\ee 
\begin{definition} 
A solitary manifold is the set ${\cal S}:=\{S(\si):\si\in\Sigma:= 
\R^3\times V\}$. 
\end{definition} 
 
The main result of our paper is the following theorem. 
\begin{theorem}\la{main} 
Let 
 (\re{ro}) and 
(\re{W}) hold. 
Let $\beta>3/2$ and $Y(t)$ be the solution to the 
Cauchy problem  (\re{WP2.1}) with 
the initial state $Y_0$ which is sufficiently close to the solitary manifold: 
\be\la{close} 
d_0:={\rm dist}_{{\cal E}_\beta}(Y_0,{\cal S})\ll 1. 
\ee 
Then the asymptotics hold for $t\to\pm\infty$, 
\be\la{qq} 
\dot q(t)=v_\pm+{\cal O}(|t|^{-2}), 
~~~~q(t)=v_\pm t+a_\pm+{\cal O}(|t|^{-3/2}); 
\ee 
\be\la{S} 
(\psi(x,t), \pi(x,t))= 
(\psi_{v_\pm}(x-v_\pm t-a_\pm), \pi_{v_\pm}(x-v_\pm t-a_\pm)) 
+W_0(t)\bPsi_\pm+r_\pm(x,t) 
\ee 
with 
\be\la{rm} 
\Vert r_\pm(t)\Vert_\cF=\cO(|t|^{-1/2}). 
\ee 
\end{theorem} 
It suffices to  prove the asymptotics  (\re{S}), (\re{qq}) for $t\to+\infty$ 
since system (\re{system0}) is time reversible.

 
\setcounter{equation}{0} 
 
\section{Symplectic Projection}

\subsection{Symplectic Structure and Hamilton Form} 
 
The system (\re{system0}) reads as the Hamilton system 
\be\la{ham} 
\dot Y=J{\cal D}{\cal H}(Y),\,\,\,J:=\left( 
\ba{cccc} 
0 & 1 & 0 & 0\\ 
-1 & 0 & 0 & 0\\ 
0 & 0 & 0 & 1\\ 
0 & 0 & -1 & 0\\ 
\ea 
\right),\,\,Y=(\psi,\pi,q,p)\in{\cal E}, 
\ee 
where ${\cal D}{\cal H}$ is the Fr\'echet derivative of 
 the Hamilton functional (\re{hamilq0}). 
Let us identify the tangent space of ${\cal E}$, at every point, 
 with the space ${\cal E}$. Consider the symplectic form $\Om$ defined on 
${\cal E}$ by the rule  
$$ 
\Om=\ds\int d\psi(x)\we d\pi(x)\,dx+dq\we dp. 
$$ 
In other words, 
\be\la{OmJ} 
\Om(Y_1,Y_2)=\langle Y_1,JY_2\rangle,\,\,\,Y_1,Y_2\in {\cal E}, 
\ee 
where 
$$ 
\langle Y_1,Y_2\rangle:=\langle\psi_1,\psi_2\rangle+ 
\langle\pi_1,\pi_2\rangle+q_1  q_2+p_1  p_2 
$$ 
and $\langle\psi_1,\psi_2\rangle=\ds\int\psi_1(x)\psi_2(x)dx$ etc. 
It is clear that the form $\Om$ is non-degenerate, i.e. 
$$ 
\Om(Y_1,Y_2)=0\,\,~\mbox{\rm for every}~~ \,Y_2\in\cE \,\,\Longrightarrow\,\, Y_1=0. 
$$ 
\begin{definition} 
i) The symbol $Y_1\nmid Y_2$ means that $Y_1\in{\cal E}$, 
$Y_2\in{\cal E}$, 
and $Y_1$ is symplectic orthogonal to $Y_2$, i.e. $\Om(Y_1,Y_2)=0$. 
 
ii) A projection operator $\bP:{\cal E}\to{\cal E}$ 
is said to be symplectic 
 orthogonal if $Y_1\nmid Y_2$ for $Y_1\in\mbox{\rm Ker}\5\bP$ and 
$Y_2\in\mbox \Im\bP$. 
\end{definition} 
 
 
\subsection{Symplectic Projection onto Solitary Manifold}

Let us consider the tangent space $\cT_{S(\si)}{\cal S}$ of 
the manifold ${\cal S}$ at a point $S(\si)$. 
The vectors $\tau_j:=\pa_{\si_j}S(\si)$, where $\pa_{\si_j}:= 
\pa_{b_j}$ and $\pa_{\si_{j+3}}:=\pa_{v_{j}}$ with $j=1,2,3$, 
form a basis in $\cT_{\si}{\cal S}$. In detail, 
\be\la{inb} 
\left.\ba{rclrrrrcrcl} 
\tau_j=\tau_j(v)&:=&\pa_{b_j}S(\si)= 
(&\!\!\!\!-\pa_j\psi_v(y)&\!\!\!\!,&\!\!\!\!-\pa_j\pi_v(y)&\!\!\!\!, 
&\!\!e_j&\!\!\!\!,&\!\!0&\!\!\!\!) 
\\ 
\tau_{j+3}=\tau_{j+3}(v)&:=&\pa_{v_j}S(\si)=(&\!\!\!\!\pa_{v_j}\psi_v(y) 
&\!\!\!\!,&\!\!\!\! 
\pa_{v_j}\pi_v(y)&\!\!\!\!,&\!\!0&\!\!\!\!,&\!\! 
\pa_{v_j}p_v&\!\!\!\!) 
\ea\right|~~~j=1,2,3, 
\ee 
where  $y:=x-b$ is the ``moving frame coordinate'', $e_1=(1,0,0)$ etc. 
Let us stress that the functions $\tau_j$  
are always regarded as  
functions of $y$ rather than those of $x$.

Formulas (\re{sol}) and conditions (\re{ro}) imply that 
\be\la{tana} 
\tau_j(v)\in\cE_\al, ~~~~v\in V, ~~~~j=1,\dots,6,~~~~~\forall\al\in\R. 
\ee 
 
\begin{lemma}\la{Ome} 
The matrix with the elements $\Om(\tau_l(v),\tau_j(v))$ is non-degenerate 
for any $v\in V$. 
\end{lemma} 
{\bf Proof } The elements are computed in Appendix A. 
As the result, the matrix $\Om(\tau_l,\tau_j)$ has the form 
\be\la{Omega} 
\Om(v):=(\Om(\tau_l,\tau_j))_{l,j=1,\dots,6}=\left( 
\ba{ll} 
0 & \Om^+(v)\\ 
-\Om^+(v) & 0 
\ea 
\right), 
\ee 
where the $3\times3$-matrix $\Om^+(v)$ equals 
\be\la{Wm} 
\Om^+(v)=K+(1-v^2)^{-1/2}E+(1-v^2)^{-3/2}v\otimes v. 
\ee 
Here $K$ is a symmetric $3\times3$-matrix with the elements 
\be\la{alpha} 
K_{ij}=\int dk|\hat\psi_v(k)|^2k_ik_j\fr{k^2+m^2+3(kv)^2}{k^2+m^2-(kv)^2} 
=\int dk|\hat\rho(k)|^2k_ik_j\fr{k^2+m^2+3(kv)^2}{(k^2+m^2-(kv)^2)^3}, 
\ee 
where the ``hat'' stands for the Fourier transform (cf. (\re{W})). 
The matrix $K$ is the integral of the symmetric 
nonnegative definite matrix 
$k\otimes 
 k=(k_ik_j)$ with a positive weight. Hence, the matrix $K$ is also 
nonnegative definite. 
Since the identity matrix $E$  
is positive definite  and the matrix $v\otimes v$ 
is nonnegative definite, the matrix $\Om^+(v)$ is  symmetric and 
positive definite, 
hence non-degenerate. Then 
the matrix $\Om(\tau_l,\tau_j)$ is also non-degenerate. 
\hfill $\bo$ 
 
Now we show that in a small neighborhood of the soliton 
manifold ${\cal S}$ a ``symplectic orthogonal projection'' 
onto ${\cal S}$ is well-defined. Let us introduce the translations 
 $T_a:(\psi(\cdot),\pi(\cdot),q,p)\mapsto 
(\psi(\cdot-a),\pi(\cdot-a),q+a,p)$, $a\in\R^3$. 
Note that the manifold ${\cal S}$ is invariant with respect to the 
translations. 
Let us denote by $v(p):=p/\sqrt{1+p^2}$ for $p\in\R^3$. 
 
\begin{definition} 
i) For any $\al\in\R$ and $\ov v<1$ denote by $\cE_\al(\ov v)=\{ 
Y=(\psi,\pi,q,p)\in\cE_\al:|v(p)|\le\ov v 
\}$. We set ${\cal E}(\ov v):={\cal E}_0(\ov v)$. 
\\ 
ii) For any  $\ti v<1$ denote by $\Si(\ti v)=\{\si=(b,v): 
b\in\R^3, |v|\le \ti v 
\}$. 
\end{definition} 
 
\begin{lemma}\la{skewpro} 
Let (\re{ro}) hold, $\al\in\R$ and $\ov v<1$. 
Then 
\\ 
i) 
there exists a neighborhood ${\cal O}_\al({\cal S})$ of ${\cal S}$ in 
${\cal E}_\al$ and a 
mapping $\bPi:{\cal O}_\al({\cal S})\to{\cal S}$  such that $\bPi$ is 
uniformly 
continuous on ${\cal O}_\al({\cal S})\cap \cE_\al(\ov v)$ 
in the metric of ${\cal E}_\al$, 
\be\la{proj} 
\bPi Y=Y~~\mbox{for}~~ Y\in{\cal S}, ~~~~~\mbox{and}~~~~~ 
Y-S \nmid \cT_S{\cal S},~~\mbox{where}~~S=\bPi Y. 
\ee 
ii) ${\cal O}_\al({\cal S})$ is invariant with respect to the translations 
 $T_a$, and 
\be\la{commut} 
\bPi T_aY=T_a\bPi Y,~~~~~\mbox{for}~~Y\in{\cal O}_\al({\cal S}) 
~~\mbox{and}~~a\in\R^3. 
\ee 
iii) For any $\ov v<1$ there exists a $\ti v<1$ s.t. 
 $\bPi Y=S(\si)$ with  $\si\in \Si(\ti v)$ for 
$Y\in{\cal O}_\al({\cal S})\cap\cE_\al(\ov v)$. 
\\\\ 
iv) For any $\ti v<1$ there exists an 
$r_\al(\ti v)>0$ s.t. 
$S(\si)+Z\in\cO_\al(\cS)$ if $\si\in\Si(\ti v)$ 
and 
$\Vert Z\Vert_\al<r_\al(\ti v)$. 
\end{lemma} 
{\bf Proof } We have to find $\si=\si(Y)$ such that $S(\si)=\bPi Y$ and 
\be\la{ift} 
\Om(Y-S(\si),\pa_{\si_j}S(\si))=0,~~~~~~ j=1,\dots,6. 
\ee 
 Let us fix an 
arbitrary $\si^0\in\Sigma$ and note that 
the system  (\re{ift}) involves only $6$ smooth scalar functions 
of $Y$. 
Then 
for $Y$ close to $S(\si^0)$, the existence of $\si$ 
follows by the standard finite dimensional 
implicit function theorem if we 
show that 
the $6\times6$ Jacobian matrix with elements $M_{lj}(Y)= 
\pa_{\si_l}\Om(Y-S(\si^0),\pa_{\si_j}S(\si^0))$ is non-degenerate 
at $Y=S(\si^0)$. First note that all the derivatives exist by  (\re{tana}). 
The non-degeneracy holds by Lemma \re{Ome} and the definition  (\re{inb}) 
since 
$M_{lj}(S(\si^0))=-\Om(\pa_{\si_l}S(\si^0),\pa_{\si_j}S(\si^0))$. 
Thus, there exists some neighborhood ${\cal O}_\al(S(\si^0))$ of $S(\si^0)$ 
where $\bPi$ is well defined and satisfies (\re{proj}), and the same is true 
in the union ${\cal O}'_\al({\cal S})=\cup_{\si^0\in \Si} 
{\cal O}_\al(S(\si^0))$. The identity (\re{commut}) holds for $Y,T_aY\in{\cal O}'_\al({\cal S})$, since the form $\Om$ and 
the manifold $\cS$ are 
invariant with respect to the translations. 
 
It remains to modify ${\cal O}'_\al({\cal S})$ by the translations: we set ${\cal O}_\al({\cal S})=\cup_{b\in\R^3}T_b{\cal O}'_\al({\cal S})$. Then the second statement obviously 
 holds. 
 
The last two statements and the uniform continuity 
in the first statement 
follow by translation invariance and 
compactness arguments. 
\hfill $\bo$ 
 
\medskip 
 
\noindent 
We refer to $\bPi$ as  
symplectic orthogonal projection onto ${\cal S}$. 
 
\bc 
The condition (\re{close}) 
implies that $Y_0=S+Z_0$ where $S=S(\si_0)=\bPi Y_0$, and 
\be\la{closeZ} 
\Vert Z_0\Vert_\beta \ll 1. 
\ee 
\ec 
\pru 
Lemma \re{skewpro} implies that $\bPi Y_0=S$ is well 
defined for small $d_0>0$. Furthermore, the condition 
(\re{close}) means that there exists a point $S_1\in\cS$ 
such that $\Vert Y_0-S_1\Vert_\beta=d_0$. 
Hence, $Y_0,S_1\in {\cal O}_\beta({\cal S})\cap\cE_\beta(\ov v)$  
with some 
$\ov v<1$ which does not depend on $d_0$ for 
sufficiently small $d_0$. 
On the other hand, $\bPi S_1= S_1$, hence 
the uniform continuity of the mapping $\bPi$ 
implies that $\Vert S_1- S\Vert_\beta\to 0$ as $d_0\to 0$. 
Therefore, finally, 
$\Vert Z_0\Vert_\beta=\Vert Y_0- S \Vert_\beta 
\le \Vert Y_0- S_1 \Vert_\beta+ 
\Vert S_1-S  \Vert_\beta\le d_0+o(1)\ll 1$ for small $d_0$.\bo 
 
 
\setcounter{equation}{0} 
 
\section{Linearization on the Solitary Manifold} 
 
Let us consider a solution to the system (\re{system0}), and split it as 
 the sum 
 
\be\la{dec} 
Y(t)=S(\si(t))+Z(t), 
\ee 
where $\si(t)=(b(t),v(t))\in\Sigma$ is an arbitrary smooth function of 
 $t\in\R$. 
In detail, denote $Y=(\psi,\pi,q,p)$ and $Z=(\Psi,\Pi,Q,P)$. 
Then (\re{dec}) means that 
\be \la{add} 
\left. 
\ba{rclrcl} 
\psi(x,t)&=&\psi_{v(t)}(x-b(t))+\Psi(x-b(t),t), 
&q(t)&=&b(t)+Q(t)\\ 
\pi(x,t)&=&\5\pi_{v(t)}(x-b(t))+\Pi(x-b(t),t), 
&p(t)&=&p_{v(t)}+P(t) 
\ea 
\right| 
\ee 
Let us 
substitute (\re{add}) to (\re{system0}), and 
linearize the equations in $Z$. Below we 
 shall choose $S(\si(t))=\bPi Y(t)$, i.e. $Z(t)$ is symplectic 
orthogonal to $\cT_{S(\si(t))}{\cal S}$. However, this orthogonality 
condition is not needed for the formal process of linearization. 
The orthogonality condition will be important in Section 6, where we derive 
``modulation equations'' for the parameters $\si(t)$.

Let us proceed to linearization. Setting $y=x-b(t)$ which is 
the ``moving frame coordinate'', we obtain from 
(\re{add}) and (\re{system0}) that 
\be\la{addeq} 
\left. 
\ba{rcl} 
\dot\psi&=& 
\dot v\cdot \na_v\psi_v(y)-\dot b\cdot \na\psi_v(y)+ 
\dot\Psi(y,t)-\dot b\cdot \na\Psi(y,t)=\pi_v(y)+\Pi(y,t) 
\\\\ 
\dot\pi&=&\dot v\cdot \na_v\pi_v(y)- 
\dot b \cdot\na\pi_v(y)+ 
\dot\Pi(y,t)-\dot b\cdot \na\Pi(y,t) 
\\\\ 
&=&\De\psi_v(y)-m^2\psi_v(y)+\De\Psi(y,t)-m^2\Psi(y,t)-\rho(y-Q) 
\\\\ 
\dot q&=&\dot b+\dot Q=\ds\fr{p_v+P}{\sqrt{1+(p_v+P)^2}} 
\\\\ 
\dot p&=&\dot v\cdot \na_v p_v+\dot P 
=-\langle\na(\psi_v(y)+ 
\Psi(y,t)),\rho(y-Q)\rangle 
\ea 
\right| 
\ee 
The equations are linear in $\Psi$ and $\Pi$, hence it remains  
extract the terms linear in $Q$ and $P$. First note that 
$\rho(y-Q)=\rho(y)-Q \cdot\na\rho(y)-N_2(Q)$, 
where 
$-N_2(Q)=\rho(y-Q)-\rho(y)+Q\cdot \na\rho(y)$. 
The condition (\re{ro}) 
implies that 
for $N_2(Q)$ the bound holds, 
\be\la{N2} 
\Vert N_2(Q)\Vert_{0,\beta}\le C_\beta(\ov Q)Q^2, 
\ee 
uniformly in $|Q|\le\ov Q$ for any fixed $\ov Q$, 
where $\beta$ is the parameter in Theorem \re{main}. 
Second, the Taylor expansion gives 
$$ 
\fr{p_v+P}{\sqrt{1+(p_v+P)^2}}=v+\nu(P-v(v \cdot P))+N_3(v,P), 
$$ 
where $\nu:=(1+p_v^2)^{-1/2}=\sqrt{1-v^2}$, and 
\be\la{N3} 
|N_3(v,P)|\le C(\ti v)P^2 
\ee 
 uniformly with respect to $|v|\le\ti v<1$. 
Using the equations (\re{stfch}),  we obtain from (\re{addeq}) 
the following equations for the components of the vector $Z(t)$: 
\be\la{Phi} 
\left.\ba{rcl} 
\!\!\!\!\dot \Psi(y,t)\!\!\!\!&=&\!\!\!\!\Pi(y,t)+\dot b\cdot \na\Psi(y,t)+ 
(\dot b-v)\cdot \na\psi_v(y)-\dot v \cdot\na_v\psi_v(y) 
\\\\ 
\!\!\!\!\dot \Pi(y,t)\!\!\!\!&=&\!\!\!\!\De\Psi(y,t)\!-\!m^2\Psi(y,t)\!+\! 
\dot b\cdot \!\na\Pi(y,t)\!+\! 
Q \!\cdot\!\na\rho(y)\!+\!(\dot b\!-\!v)\!\cdot\! \na\pi_v(y)\!- 
\!\dot v\cdot \!\na_v\pi_v(y)\!+\!N_2 
\\\\ 
\!\!\!\!\dot Q(t)\!\!\!\!&=&\!\!\!\!\nu (E-v\otimes v)P+(v-\dot b)+N_3 
\\\\ 
\!\!\!\!\dot P(t)\!\!\!\!&=&\!\!\!\!\langle\Psi(y,t),\na\rho(y)\rangle+ 
\langle\na\psi_v(y), 
Q \cdot\na\rho(y)\rangle-\dot v\cdot \na_v p_v+N_4(v,Z) 
\ea\right| 
\ee 
where $N_4(v,Z)=\langle\na\psi_v,N_2(Q)\rangle+ 
\langle\na\Psi,Q\cdot \na\rho\rangle+\langle\na\Psi,N_2(Q)\rangle$. Clearly, 
$N_4(v,Z)$ satisfies the following estimate 
\be\la{N4} 
|N_4(v,Z)|\le C_\beta(\rho,\ti v,\ov Q)\Big[Q^2+ 
\Vert\Psi\Vert_{1,-\beta}|Q| 
 \Big], 
\ee 
uniformly in $|v|\le \ti v$ and $|Q|\le \ov Q$ 
for any fixed  $\ti v<1$. 
We can write the equations (\re{Phi}) as 
\be\la{lin} 
\dot Z(t)=A(t)Z(t)+T(t)+N(t),\,\,\,t\in\R. 
\ee 
Here the operator $A(t)=A_{v,w}$ depends on two parameters, $v=v(t)$, and 
$w=\dot b(t)$ and can be written in the form 
\be\la{AA} 
A_{v,w}\left( 
\ba{c} 
\Psi \\ \Pi \\ Q \\ P 
\ea 
\right):=\left( 
\ba{cccc} 
w \cdot\na & 1 & 0 & 0 \\ 
\De-m^2 & w\cdot \na & \na\rho\cdot & 0 \\ 
0 & 0 & 0 & B_v \\ 
\langle\cdot,\na\rho\rangle & 0 & \langle\na\psi_v,\cdot\na\rho\rangle & 0 
\ea 
\right)\left( 
\ba{c} 
\Psi \\ \Pi \\ Q \\ P 
\ea 
\right), 
\ee 
where $B_{v}=\nu(E-v\otimes v)$. 
Furthermore,   $T(t)=T_{v,w}$ and $N(t)=N(\si,Z)$ 
are given by 
\be\la{TN} 
T_{v,w}=\left( 
\ba{c} 
(w-v)\cdot\na\psi_v-\dot v\cdot\na_v\psi_v\\ 
(w-v)\cdot\na\pi_v-\dot v\cdot\na_v\pi_v\\ 
v-w \\ 
-\dot v\cdot\na_v p_v 
\ea 
\right),\,\,\,\,N(\si,Z)=\left( 
\ba{c} 
0 \\ N_2(Z) \\ N_3(v,Z) \\ N_4 (v,Z) 
\ea 
\right), 
\ee 
where $v=v(t)$, $w=w(t)$, $\si=\si(t)=(b(t),v(t))$, and $Z=Z(t)$. 
Estimates (\re{N2}), (\re{N3}) and (\re{N4}) imply that 
\be\la{N14} 
\Vert N(\si,Z)\Vert_\beta\le C(\ti v, \ov Q) 
\Vert Z\Vert_{-\beta}^2, 
\ee 
uniformly in $\si\in\Si(\ti v)$ and $\Vert Z\Vert_{-\beta}\le r_{-\beta}(\ti v)$ 
for any fixed  $\ti v<1$. 
\brs\la{rT} 
{\rm 
i) 
The term 
$A(t)Z(t)$ 
in the right hand side of equation  (\re{lin}) 
is linear  in $Z(t)$, 
and $N(t)$ is a {\it high order term} in $Z(t)$. 
On the other hand, 
$T(t)$ is a zero order term which does not vanish at $Z(t)=0$ 
since 
$S(\si(t))$ generally is not a soliton solution if (\re{sigma}) 
fails to hold (though $S(\si(t))$ belongs to the solitary manifold). 
\\ 
ii) Formulas (\re{inb}) and (\re{TN}) imply: 
\be\la{Ttang} 
T(t)=-\sum\limits_{l=1}^3[(w-v)_l\tau_l+\dot v_l\tau_{l+3}] 
\ee 
and hence $T(t)\in \cT_{S(\si(t))}{\cal S}$, $t\in\R$. 
This fact suggests an unstable character of the nonlinear dynamics 
{\it along the solitary manifold}. 
} 
\ers 
 
 \setcounter{equation}{0} 
\section{Linearized Equation} 
Here we collect some Hamiltonian and spectral properties of the 
operator (\re{AA}). 
 First, let us consider the linear equation 
\be\la{line} \dot X(t)=A_{v,w}X(t),~~~~~~~t\in\R \ee with arbitrary 
fixed $v\in V=\{v\in\R^3: |v|<1\}$ and $w\in \R^3$. Let us define 
the space $\cE^+:=H^2(\R^3)\oplus H^1(\R^3)\oplus\R^3\oplus\R^3$.

\begin{lemma} \la{haml} 
i) For any $v\in V$ and $w\in \R^3$ 
equation (\re{line}) 
 can be represented as the Hamiltonian system  (cf. (\re{ham})), 
\be\la{lineh} 
\dot X(t)= 
JD{\cal H}_{v,w}(X(t)),~~~~~~~t\in\R, 
\ee 
where $D{\cal H}_{v,w}$ is the Fr\'echet derivative of the 
Hamiltonian functional 
\beqn 
{\cal H}_{v,w}(X)&=&\fr12\int\Big[|\Pi|^2+ 
|\na\Psi|^2+m^2|\Psi|^2\Big]dy+\int\Pi w\cdot\na\Psi dy+\int\rho(y) 
 Q\cdot\na\Psi dy 
\nonumber\\ 
&+& 
\fr12P\cdot B_vP-\fr12\langle Q\cdot\na\psi_{v}(y),Q\cdot\na\rho(y)\rangle, 
~~~~X=(\Psi,\Pi,Q,P)\in \cE. 
\la{H0} 
\eeqn 
\\ 
ii) The energy conservation law holds 
for the solutions $X(t)\in C^1(\R,\cE^+)$, 
\be\la{enec} 
\cH_{v,w}(X(t))=\co,~~~~~t\in\R. 
\ee 
iii) The skew-symmetry relation holds, 
\be\la{com} 
\Omega(A_{v,w}X_1,X_2)=-\Omega(X_1,A_{v,w}X_2), ~~~~~~~~X_1,X_2\in \cE. 
\ee 
\end{lemma} 
{\bf Proof } i) The equation (\re{line}) reads as follows, 
\be\la{eql} 
\fr{d}{dt}\left( 
\ba{c} 
\Psi \\ \Pi \\ Q \\ P 
\ea 
\right)=\left( 
\ba{l} 
\Pi+w\cdot\na\Psi \\ 
\De\Psi-m^2\Psi+w\cdot\na\Pi+Q\cdot\na\rho \\ 
B_{v}P \\ 
-\langle\na\Psi,\rho\rangle+\langle\na\psi_{v},Q\cdot\na\rho\rangle 
\ea 
\right). 
\ee 
The first three equations correspond to the Hamilton form 
since 
$$ 
\Pi+w\cdot\na\Psi=D_\Pi{\cal H}_{v,w},\,\,\, \De\Psi- 
m^2\Psi+w\cdot\na\Pi+Q\cdot\na\rho=-D_\Psi{\cal H}_{v,w}, 
\,\,\, 
B_{v}P=\na_P{\cal H}_{v,w}. 
$$ 
Let us check that the last equation 
has also the Hamilton form, i.e. 
$-\langle\na\Psi,\rho\rangle+ 
\langle\na\psi_{v},Q\cdot\na\rho\rangle=-\na_Q{\cal H}_{v,w}$. 
First we note that $-\langle\pa_j\Psi,\rho\rangle=-\pa_{Q_j} 
\ds\int \rho Q\cdot\na\Psi dx$. It remains to show that 
\be\la{paq} 
\langle\pa_j\psi_{v},Q\cdot\na\rho\rangle=\pa_{Q_j} 
\fr12\langle Q\cdot\na\psi_{v},Q\cdot\na\rho\rangle. 
\ee 
Indeed, 
\beqn 
\pa_{Q_j}(\fr12\langle Q\cdot\na\psi_{v},Q\cdot\na\rho\rangle&=& 
\fr12\langle\pa_j\psi_v,Q\cdot\na\rho\rangle+\fr12\langle Q\cdot\na\psi_v, 
\pa_j\rho\rangle 
\nonumber\\ 
&=& \fr12\langle\pa_j\psi_v,Q\cdot\na\rho\rangle+\fr12\langle\pa_j\psi_{v}, 
Q\cdot\na\rho\rangle 
\eeqn 
where we have 
integrated  twice by parts. 
Then (\re{paq}) follows. 
\\ 
ii) The energy conservation law 
follows by  (\re{lineh}) and the chain rule for the 
Fr\'echet derivatives: 
\be\la{crF} 
\ds\fr d{dt}\cH_{v,w} 
(X(t))=\langle D\cH_{v,w}(X(t)),\dot X(t)\rangle= 
\langle D\cH_{v,w}(X(t)),J D\cH_{v,w}(X(t))\rangle=0,~~~~~~t\in\R 
\ee 
since the operator $J$ is skew-symmetric by (\re{ham}), and 
$D\cH_{v,w}(X(t))\in \cE$ for  $X(t)\in \cE^+$. 
\\ 
iii) The skew-symmetry 
holds 
since $A_{v,w}X=JD\cH_{v,w}(X)$, 
and the linear operator 
$X\mapsto D\cH_{v,w}(X)$ 
is symmetric as the Fr\'echet derivative of a quadratic form. 
\hfill$\bo$ 
 
\br 
{\rm 
One can obtain (\re{H0}) by  expanding ${\cal H}(S_{b,v}+X)$ to a power 
 series in $X$ up to second order terms. As a result, ${\cal H}_{v,w}(X)$ 
is the quadratic part of the Taylor series complemented by 
the second integral on the right hand side of (\re{H0}) arising from the 
left hand side of (\re{ham}). 
} 
\er 
\begin{lemma} \la{ljf} 
The operator $A_{v,w}$ acts on the 
tangent vectors $\tau_j(v)$ to the solitary manifold 
as follows, 
\be\la{Atanform} 
A_{v,w}[\tau_j(v)]=(w-v)\cdot\na\tau_j(v),\,\,\,A_{v,w}[\tau_{j+3}(v)]= 
(w-v)\cdot\na\tau_{j+3}(v)+\tau_j(v),\,\,\,j=1,2,3. 
\ee 
\end{lemma} 
{\bf Proof } In detail, we have to show that 
$$ 
A_{v,w}\left( 
\ba{c} 
-\pa_j\psi_{v} \\ -\pa_j\pi_{v}  \\ e_j \\ 0 
\ea 
\right)=\left( 
\ba{c} 
(v-w)\cdot\na\pa_j\psi_{v}  \\ (v-w)\cdot\na\pa_j\pi_{v}  \\ 0 \\ 0 
\ea 
\right), 
$$ 
\be\la{Atan} 
A_{v,w}\left( 
\ba{c} 
\pa_{v_j}\psi_{v}  \\ \pa_{v_j}\pi_{v}  \\ 0 \\ \pa_{v_j}p_{v} 
\ea 
\right)=\left( 
\ba{c} 
(w-v)\cdot\na\pa_{v_j}\psi_{v}  \\ (w-v)\cdot\na\pa_{v_j}\pi_{v} 
\\ 0 \\ 0 
\ea 
\right)+\left( 
\ba{c} 
-\pa_{j}\psi_{v}  \\ -\pa_{j}\pi_{v}  \\ e_j \\ 0 
\ea 
\right). 
\ee 
Indeed, differentiate the equations (\re{stfch}) in $b_j$ 
 and $v_j$, and obtain that the derivatives of soliton state in 
parameters satisfy the following equations, 
\be\la{stinb} 
\left.\ba{rclrcl} 
-v\cdot\na\pa_j\psi_{v}\!\!\!\!&=&\!\!\!\!\pa_j\pi_{v}\,, 
& 
-v\cdot\na\pa_j\pi_{v}\!\!\!\!&=&\!\!\!\!\De\pa_j\psi_{v}-m^2\pa_j\psi_{v}-\pa_j\rho 
\\\\ 
-\pa_j\psi_{v}-v\cdot\na\pa_{v_j}\psi_{v}\!\!\!\!&=&\!\!\!\!\pa_{v_j}\pi_{v}\,, 
& 
-\pa_j\pi_{v}- v\cdot\na\pa_{v_j}\pi_{v}\!\!\!\!&=&\!\!\!\! 
\De\pa_{v_j}\psi_{v}-m^2\pa_{v_j}\psi_{v}\,\\\\ 
\pa_{v_j}p_{v}= e_j(1-v^2)^{-1/2}\!\!\!\!&+&\!\!\!\!\ds v\fr{v_j}{(1-v^2)^{3/2}}\,,& 
0\!\!\!\!&=&\!\!\!\!-\langle\na\pa_{v_j}\psi_{v},\rho\rangle 
\ea\right| 
\ee 
for $j=1,2,3$. 
Then (\re{Atan}) follows from (\re{stinb}) by definition 
of $A$ in (\re{AA}).\hfill$\bo$

We shall apply Lemma \re{haml} mainly to the operator $A_{v,v}$ 
corresponding to $w=v$. 
In that case the linearized equation has the following 
additional specific 
features. 
 
 \bl\la{ceig} 
Let us assume that $w=v\in V$. Then 
\\ 
i) 
The tangent vectors $\tau_j(v)$ with $j=1,2,3$ are eigenvectors, 
and $\tau_{j+3}(v)$ are root vectors of the 
operator $A_{v,v}$, that correspond to the zero eigenvalue, i.e. 
\be\la{Atanformv} 
A_{v,v}[\tau_j(v)]=0,\,\,\,A_{v,v}[\tau_{j+3}(v)]= 
\tau_j(v),\,\,\,j=1,2,3. 
\ee 
ii) 
The Hamiltonian function (\re{H0}) is  nonnegative definite since 
\be\la{H0vv} 
{\cal H}_{v,v}(X)=\ds\fr12\int\Big[|\Pi+v\cdot\na\Psi|^2+ 
|\Lam^{1/2}\Psi-\Lam^{-1/2}Q\cdot\na\rho|^2\Big]dx+ 
\fr12P\cdot B_vP\ge 0. 
\ee 
Here $\Lam$ is the operator (\re{Lpv}) which is symmetric and nonnegative 
definite in $L^2(\R^3)$ 
for $|v|<1$, 
and $\Lam^{1/2}$ is the nonnegative definite square root defined in the Fourier 
representation. 
\el 
\pru 
The 
first statement follows from 
(\re{Atanform}) with $w=v$. 
In order to prove ii) we 
rewrite the integral in (\re{H0vv})  as follows: 
\beqn\la{H0vve} 
&&\ds\fr 12\langle \Pi+v\cdot\na\Psi,\Pi+v\cdot\na\Psi \rangle+ 
\fr 12 
\langle 
\Lam^{1/2}\Psi-\Lam^{-1/2}Q\cdot\na\rho, 
\Lam^{1/2}\Psi-\Lam^{-1/2}Q\cdot\na\rho\rangle 
\nonumber\\\nonumber\\ 
&=&\ds\fr 12 
\langle \Pi,\Pi \rangle+ 
\langle \Pi,v\cdot\na\Psi \rangle+ 
\ds\fr 12 
\langle v\cdot\na\Psi,v\cdot\na\Psi \rangle 
\nonumber\\\nonumber\\ 
&+& 
\ds\fr 12 
\langle 
\Lam\Psi, 
\Psi\rangle- 
\langle 
\Psi, 
Q\cdot\na\rho\rangle+ 
\ds\fr 12 
\langle 
\Lam^{-1}Q\cdot\na\rho, 
Q\cdot\na\rho\rangle 
\eeqn 
since the operator $\Lam^{1/2}$ is symmetric in $L^2(\R^3)$. 
Now all the terms of 
the expression (\re{H0vve}) can be identified with the 
corresponding terms in (\re{H0}) since 
\beqn\la{H0vvs} 
&&\ds\fr 12 
\langle 
\Lam\Psi, 
\Psi\rangle= 
\ds\fr 12 
\langle 
[-\De+m^2+(v\cdot\na)^2]\Psi, 
\Psi\rangle,~~~~~~~~~~ 
\Lam^{-1}\rho=-\psi_v 
\eeqn 
by (\re{Lpv}) and (\re{Lpvy}).\bo 
\br 
{\rm 
In Section 14 
we will apply Lemma \re{ceig} ii) together with energy conservation 
(\re{enec}) to prove the analyticity of the resolvent $(A_{v,v}-\lam)^{-1}$ 
for $\Re\lam>0$. 
} 
\er

\br 
{\rm 
 For a soliton solution of the system 
(\re{system0}) we have $\dot b=v$, $\dot v=0$, and hence $T(t)\equiv 0$. 
Thus, equation 
(\re{line}) 
is the linearization of system 
(\re{system0}) on a soliton solution.  
In fact, we linearize (\re{system0}) on a trajectory 
$S(\si(t))$, where $\si(t)$ is nonlinear with respect to $t$, 
rather than on a soliton solution.  
We shall show below  
that $T(t)$ is quadratic in $Z(t)$ 
if we choose 
$S(\si(t))$ to be  the symplectic orthogonal projection of $Y(t)$. 
In this case,  
(\re{line})  is a linearization of (\re{system0}) again. 
} 
\er 
 
 
\setcounter{equation}{0} 
 
\section{Symplectic Decomposition of the Dynamics} 
 
Here we decompose the dynamics in two components: along the manifold 
 ${\cal S}$ and in transversal directions. The equation (\re{lin}) 
is obtained without any assumption on $\si(t)$ in (\re{dec}). 
We are going to specify $S(\si(t)):=\bPi Y(t)$. 
However, in this case we must known that 
\be\la{YtO} 
Y(t)\in \cO_\al(\cS),~~~~~t\in\R, 
\ee 
with some $\cO_\al(\cS)$ defined in Lemma \re{skewpro}. 
It is true for $t=0$ by our main assumption 
 (\re{close}) with sufficiently small $d_0>0$. 
Then  $S(\si(0))=\bPi Y(0)$ and  $Z(0)=Y(0)-S(\si(0))$ 
are well defined. 
We shall prove below that (\re{YtO}) holds 
with $\al=-\beta$ 
 if $d_0$ is sufficiently small. 
First, the a priori estimate (\re{WP2.1'}) 
together with Lemma \re{skewpro} iii) 
imply that $\bPi Y(t)=S(\si(t))$ with $\si(t)=(b(t),v(t))$, 
 and 
\be\la{vsigmat} 
|v(t)|\le \ti v<1,~~~~~~t\in\R 
\ee 
if $Y(t)\in \cO_{-\beta}(\cS)$. Denote by $r_{-\beta}(\ti v)$  
the positive 
number in 
Lemma \re{skewpro} iv) which corresponds to $\al=-\beta$. 
Then $S(\si)+Z\in \cO_{-\beta}(\cS)$ if 
$\si=(b,v)$ with $|v|<\ti v$ and 
$ \Vert Z\Vert_{-\beta}<r_{-\beta}(\ti v)$. 
Note that (\re{WP2.1'}) implies 
$\Vert Z(0)\Vert_{-\beta}<r_{-\beta}(\ti v)$ if 
 $d_0$ is sufficiently small. 
Therefore, $S(\si(t))=\bPi Y(t)$ and  $Z(t)=Y(t)-S(\si(t))$ 
are well defined for small times 
$t\ge 0$, 
such that 
 $\Vert Z(t)\Vert_{-\beta} < r_{-\beta}(\ti v)$. 
This  
argument can be  
 formalized by the following standard definition. 
\begin{definition} 
Let $t_*$ be the ``exit time'', 
\be\la{t*} 
t_*=\sup \{t>0: 
\Vert Z(s)\Vert_{-\beta} < r_{-\beta}(\ti v),~~0\le s\le t\}. 
\ee 
\end{definition} 
 
One of our main goals is to prove that $t_*=\infty$ 
if $d_0$ is sufficiently small. 
This would follow if we shall show that 
\be\la{Zt} 
\Vert Z(t)\Vert_{-\beta}<r_{-\beta}(\ti v)/2,~~~~~0\le t < t_*. 
\ee 
Note that 
\be\la{Qind} 
|Q(t)|\le\ov Q:= r_{-\beta}(\ti v), ~~~~~0\le t< t_*. 
\ee 
Now  by (\re{N14}), 
the term $N(t)$ in (\re{lin}) 
satisfies the following estimate, 
\be\la{Nest} 
\Vert N(t)\Vert_{\beta}\le C_\beta(\ti v)\Vert Z(t)\Vert^2_{-\beta}, 
\,\,\,0\le t<t_*. 
\ee

 
\subsection{Longitudinal Dynamics: Modulation Equations} 
 
From now on 
 we fix the decomposition 
$Y(t)=S(\si(t))+Z(t)$ for $0<t<t_*$ 
by setting $S(\si(t))=\bPi Y(t)$ which is equivalent to the 
 symplectic orthogonality condition of type (\re{proj}), 
\be\la{ortZ} 
Z(t)\nmid\cT_{S(\si(t))}{\cal S},\,\,\,0\le t<t_*. 
\ee 
This enables us drastically simplify  the 
asymptotic analysis of the 
dynamical equation (\re{lin}) 
 for the transversal component $Z(t)$.  
As the first step, we derive the 
longitudinal dynamics, i.e. find  
the ``modulation equations'' for the parameters 
 $\si(t)$. 
Let us derive a system of ordinary differential equations for the 
vector $\si(t)$. For this purpose, 
let us write (\re{ortZ}) in the form 
\be\la{orth} 
\Om(Z(t),\tau_j(t))=0,\,\,j=1,\dots,6, ~~~~~~~0\le t<t_*, 
\ee 
where the vectors $\tau_j(t)=\tau_j(\si(t))$ span the 
tangent space 
$\cT_{S(\si(t))}{\cal S}$. 
Note that $\si(t)=(b(t),v(t))$, where 
\be\la{sit} 
|v(t)|\le \ti v<1,~~~~~~~~~0\le t<t_*, 
\ee 
by Lemma \re{skewpro} iii). 
  It would be convenient for us to use some other 
parameters $(c,v)$ instead of $\si=(b,v)$, where 
$c(t)= 
b(t)-\ds\int^t_0 v(\tau)d\tau$ and 
\be\la{vw} 
\dot c(t)=\dot b(t)-v(t)=w(t)-v(t), ~~~~~~~~~0\le t<t_*. 
\ee 
We do not 
need an explicit form of the equations for $(c,v)$ but the following 
 statement. 
 
\begin{lemma}\la{mod} 
Let $Y(t)$ be a solution to the Cauchy problem (\re{WP2.1}), 
and (\re{dec}), 
(\re{orth}) hold. Then $(c(t),v(t))$ satisfies 
the equation 
\be\la{parameq} 
\left( 
\ba{l} 
\dot c(t) \\ \dot v(t) 
\ea 
\right)={\cal N}(\si(t),Z(t)), ~~~~~~~0\le t<t_*, 
\ee 
where 
\be\la{NZ} 
{\cal N}(\si,Z)={\cal O}(\Vert Z\Vert^2_{-\beta}) 
\ee 
uniformly in  $\si\in\Si(\ti v)$. 
\end{lemma} 
{\bf Proof }  We differentiate the orthogonality conditions 
 (\re{orth})  in $t$, and obtain 
\be\la{ortder} 
0=\Omega(\dot Z,\tau_j)+\Omega(Z,\dot\tau_j)=\Omega(AZ+T+N,\tau_j)+ 
\Omega(Z,\dot\tau_j), ~~~~~~~0\le t<t_*. 
\ee 
First, let us compute the principal (i.e. non-vanishing at $Z=0$) 
term $\Om(T,\tau_j)$. For 
$j=1,2,3$ one has by (\re{Ttang}), (\re{Omega}) 
$$ 
\Omega(T,\tau_j)=-\sum\limits_l(\dot c_l\Om(\tau_l,\tau_j)+ 
\dot v_l\Om(\tau_{l+3},\tau_j))= 
\sum\limits_l\Om(\tau_j,\tau_{l+3})\dot v_l=\sum\limits_l\Om^+_{jl}\dot v_l, 
$$ 
where the matrix $\Om^+$ is defined by (\re{Wm}). 
Similarly, 
$$ 
\Om(T,\tau_{j+3})=-\sum\limits_l(\dot c_l\Om(\tau_l,\tau_{j+3})+ 
\dot v_l\Om(\tau_{l+3},\tau_{j+3}))= 
\sum\limits_l\Om(\tau_{j+3},\tau_l)\dot c_l=-\sum\limits_l\Om^+_{jl}\dot c_l. 
$$ 
As the result, we have by (\re{Omega}), 
\be\la{OmTtau} 
\Om(T,\tau)= 
\left( 
\ba{ll} 
0 & \Om^+(v) \\ 
-\Om^+(v) & 0 
\ea 
\right) 
\left( 
\ba{c} 
\dot c \\ \dot v 
\ea 
\right)=\Om(v) 
\left( 
\ba{c} 
\dot c \\ \dot v 
\ea 
\right) 
\ee 
 in the vector form. 
 
Second, let us compute $\Omega(AZ,\tau_j)$. 
The skew-symmetry (\re{com}) implies that 
$\Omega(AZ,\tau_j)=-\Omega(Z,A\tau_j)$. 
Then for $j=1,2,3$, we have by (\re{Atanform}), 
\be\la{omaz1} 
\Omega(AZ,\tau_j)=-\Omega(Z,\dot c\cdot\na\tau_j), 
\ee 
and similarly, 
\be\la{omaz2} 
\Omega(AZ,\tau_{j+3})=-\Omega(Z,\dot c\cdot\na\tau_{j+3}+\tau_j)= 
-\Omega(Z,\dot c\cdot\na\tau_{j+3})-\Omega(Z,\tau_j)= 
-\Omega(Z,\dot c\cdot\na\tau_{j+3}), 
\ee 
since $\Om(Z,\tau_j)=0$. 
 
Finally, let us compute the last term $\Om(Z,\dot\tau_j)$. For 
$j=1,\dots,6$ one has $\dot\tau_j=\dot b\cdot\na_b\tau_j+ \dot 
v\cdot\na_v\tau_j= \dot v\cdot\na_v\tau_j$ since the vectors 
$\tau_j$ do not depend on $b$ according to (\re{inb}). 
Hence, 
\be\la{Ztau} 
\Omega(Z,\dot\tau_j)= \Om(Z,\dot v\cdot\na_v\tau_j). 
\ee 
As the result, by (\re{OmTtau})-(\re{Ztau}), 
 the equation (\re{ortder}) becomes 
\be\la{modul} 
0=\Om(v) 
\left( 
\ba{l} 
\dot c \\ \dot v 
\ea 
\right)+{\cal M}_0(\si,Z)\left( 
\ba{l} 
\dot c \\ \dot v 
\ea 
\right)+{\cal N}_0(\si,Z), 
\ee 
where the matrix ${\cal M}_0(\si,Z)={\cal O}(\Vert Z\Vert_{-\beta})$, 
and 
${\cal N}_0(\si,Z)={\cal O}(\Vert Z\Vert^2_{-\beta})$ 
uniformly in $\si\in\Si(\ti v)$ and 
$\Vert Z\Vert_{-\beta}<r_{-\beta}(\ti v)$. 
Then, since $\Om(v)$ is invertible 
by Lemma \re{Ome}, 
and $\Vert Z\Vert_{-\beta}$ is small, we can resolve (\re{modul}) 
 with respect to the derivatives and obtain equations (\re{parameq}) 
 with ${\cal N}={\cal O}(\Vert Z\Vert^2_{-\beta})$ uniformly 
in $\si\in\Sigma(\ti v)$.\hfill$\bo$ 
\br \la{radiab} 
{\rm 
The equations (\re{parameq}), (\re{NZ}) imply that the soliton 
parameters $c(t)$ and $v(t)$ are {\it adiabatic invariants} (see \ci{AKN}). 
} 
\er 
 
\subsection{Decay for the Transversal Dynamics} 
 
In Section 11 we shall show that our main Theorem \re{main} 
can be derived from the following time decay of the 
transversal component $Z(t)$: 
\bp\la{pdec} 
 Let all conditions of Theorem \re{main} hold. Then $t_*=\infty$, and 
\be\la{Zdec} 
\Vert Z(t)\Vert_{-\beta}\le \ds\fr {C(\rho,\ov v,d_0)}{(1+|t|)^{3/2}},~~~~~t\ge0. 
\ee 
\ep 
We shall derive (\re{Zdec}) in Sections 7-11 from our equation 
(\re{lin}) for the transversal component $Z(t)$. 
This equation can be specified by  
using Lemma \re{mod}. 
Indeed, the lemma implies that 
\be\la{Tta} 
\Vert T(t)\Vert_{\beta}\le C(\ti v)\Vert Z(t)\Vert^2_{-\beta}, 
~~~~~~~~~0\le t<t_*, 
\ee 
by (\re{TN})  since $w-v=\dot c$. 
Thus, 
equation (\re{lin}) becomes  
\be\la{reduced} 
\dot Z(t)=A(t)Z(t)+\ti N(t), ~~~~~~~~~0\le t<t_*, 
\ee 
where 
$A(t)=A_{v(t),w(t)}$, and 
$\ti N(t):=T(t)+N(t)$ satisfies the estimate 
\be\la{redN} 
\Vert\ti  N(t)\Vert_{\beta}\le C\Vert Z(t)\Vert^2_{-\beta},~~~~ 
~~~~~~~~~0\le t<t_*. 
\ee 
In the remaining part of our paper we mainly analyze the 
{\bf basic equation} (\re{reduced}) to establish the decay (\re{Zdec}). 
We are going to derive the decay 
using the bound (\re{redN}) and 
 the 
orthogonality condition  (\re{ortZ}). 
 
Let us comment on two main difficulties in proving 
(\re{Zdec}). The difficulties are common for the problems 
studied in \ci{BP1}. 
First, the linear part of the equation is non-autonomous, 
hence we cannot apply directly known methods of scattering theory. 
Similarly to the approach of  \ci{BP1}, 
we reduce the problem to 
 the analysis of the 
 {\it frozen} linear equation, 
\be\la{Avv} \dot X(t)=A_1X(t), ~~t\in\R, 
\ee 
where $A_1$ is the operator 
$A_{v_1,v_1}$ defined by 
(\re{AA}) 
with $v_1=v(t_1)$ for a fixed $t_1\in[0,t_*)$. Then we 
estimate the error by the method of majorants. 
 
Second, even for the frozen equation (\re{Avv}), 
the decay 
of type  (\re{Zdec}) for all solutions does not hold 
 without  the 
orthogonality condition  of type (\re{ortZ}). 
Namely, by  (\re{Atanformv}) 
equation 
(\re{Avv}) admits the {\it secular solutions} 
\be\la{secs} 
X(t)=\sum_1^3 C_{j}\tau_j(v_1)+\sum_1^3 D_j[\tau_j(v_1)t+ 
\tau_{j+3}(v_1)]. 
\ee 
 
The solutions lie in the tangent space   
$\cT_{S(\si_1)}{\cal S}$ 
with 
$\si_1=(b_1,v_1)$ (for an arbitrary $b_1\in\R$) 
that 
suggests an unstable character of the nonlinear dynamics 
{\it along the solitary manifold} (cf. Remark \re{rT} ii)). 
Thus, the 
orthogonality condition  (\re{ortZ}) eliminates   
the secular solutions. 
We shall apply the corresponding 
projection to kill the unstable  
``longitudinal terms'' in the basic equation 
 (\re{reduced}).

\bd 
i) For $v\in V$, denote by $\bPi_v$ the symplectic orthogonal projection 
of ${\cal E}$ onto the tangent space $\cT_{S(\si)}{\cal S}$, and 
write $\bP_v=\bI-\bPi_v$. 
\\ 
ii) Denote by $\cZ_v=\bP_v\cE$ the space symplectic 
orthogonal to $\cT_{S(\si)}{\cal S}$ with 
$\si=(b,v)$ (for an arbitrary $b\in\R$). 
\ed 
Note that by the linearity, 
\be\la{Piv} 
\bPi_vZ=\sum\bPi_{jl}(v) 
\tau_j(v)\Om(\tau_l(v),Z),~~~~~~~~~~Z\in\cE, 
\ee 
 with some smooth coefficients $\bPi_{jl}(v)$. 
Hence, the projector $\bPi_v$  
does not depend on $b$ 
(in the variable $y=x-b$), 
and this explains the choice 
of the subindex in $\bPi_v$ and $\bP_v$. 
 
We have now the symplectic orthogonal decomposition 
\be\la{sod} 
\cE=\cT_{S(\si)}{\cal S}+\cZ_v,~~~~~~~\si=(b,v), 
\ee 
and the symplectic orthogonality  (\re{ortZ}) 
can be represented in the following equivalent forms, 
\be\la{PZ} 
\bPi_{v(t)} Z(t)=0,~~~~\bP_{v(t)}Z(t)= Z(t),~~~~~~~~~0\le t<t_*. 
\ee 
 
\br\la{rZ} 
{\rm 
The tangent space $\cT_{S(\si)}{\cal S}$ is invariant under 
the operator $A_{v,v}$ by Lemma \re{ceig} i), hence 
the space  $\cZ_v$ is also invariant by 
(\re{com}), namely: $A_{v,v}Z\in \cZ_v$ 
on a dense domain of $Z\in \cZ_v$. 
} 
 
\er

In Sections 12-18 below we will prove the following 
 proposition which is one of the main ingredients to proving 
(\re{Zdec}). Let us 
consider the Cauchy problem for equation (\re{Avv}) 
with $A=A_{v,v}$ for a fixed $v\in V$. 
Recall that the parameter $\beta>3/2$ is also fixed. 
\begin{pro}\la{lindecay} 
Let  (\re{ro}) and (\re{W}) hold, 
$|v|\le\ti v<1$, and $X_0\in\cE$. 
Then 
\\ 
i) The equation (\re{Avv}), with $A_1=A=A_{v,v}$, admits a 
unique solution 
 $e^{At}X_0:=X(t)\in C(\R, \cE)$ 
with the initial condition $X(0)=X_0$. 
\\ 
ii) 
For $X_0\in 
\cZ_{v}\cap \cE_\beta$, the solution  
$X(t)$ has the 
following decay, 
\be\la{frozenest} 
\Vert e^{At}X_0\Vert_{-\beta}\le 
\fr{C(\beta,\ti v)}{(1+|t|)^{3/2}}\Vert X_0\Vert_{\beta},~~~~~~~~ 
\,\,\,t\in\R. 
\ee 
\end{pro} 
\begin{remark}{\rm 
The decay is provided by two fundamental facts 
which we will establish below: 
\\ 
i) the null root space of the generator $A$ 
coincides with the tangent space 
$\cT_{S(\si)}{\cal S}$, 
where 
$\si=(b,v)$ (for an arbitrary $b\in\R$), 
and 
\\ 
ii) the  spectrum of $A$ 
in the space $\cZ_{v}$ is purely continuous.} 
\end{remark}

\setcounter{equation}{0} 
 
\section{Frozen Form of Transversal Dynamics}

Now  let us fix an arbitrary $t_1\in [0,t_*)$, and 
rewrite the equation (\re{reduced}) in a ``frozen form'' 
\be\la{froz} 
\dot Z(t)=A_1Z(t)+(A(t)-A_1)Z(t)+\ti N(t),\,\,\,~~~~0\le t<t_*, 
\ee 
where $A_1=A_{v(t_1),v(t_1)}$ and 
$$ 
A(t)-A_1=\left( 
\ba{cccc} 
[w(t)\!-\!v(t_1)]\cdot \na & 0 & 0 & 0 \\ 
0 & [w(t)\!-\!v(t_1)]\cdot \na & 0 & 0 \\ 
0 & 0 & 0 & B_{v(t)}\!-\!B_{v(t_1)} \\ 
0 & 0 & \langle\na(\psi_{v(t)}\!-\!\psi_{v(t_1)}),\na\rho\rangle & 0 
\ea 
\right). 
$$ 
The next trick is important since it enables 
 us to kill the ``bad terms'' 
 $[w(t)\!-\!v(t_1)]\cdot \na$ in the operator $A(t)-A_1$. 
\begin{definition}\la{d71} 
Let us change the  variables $(y,t)\mapsto (y_1,t)=(y+d_1(t),t)$ 
where 
\be\la{dd1} 
d_1(t):=\int_{t_1}^t(w(s)-v(t_1))ds, ~~~~0\le t\le t_1. 
\ee 
\end{definition} 
Next, let us write 
\be\la{Z1} 
Z_1(t)= 
(\Psi(y_1-d_1(t),t),\Pi(y_1-d_1(t),t),Q(t),P(t)). 
\ee 
Then we obtain the final form of the 
``frozen equation'' for the transversal dynamics 
\be\la{redy1} 
\dot Z_1(t)=A_1Z_1(t)+B_1(t)Z_1(t)+N_1(t),\,\,\,0\le t\le t_1, 
\ee 
where $N_1(t)=\ti N(t)$ expressed in terms of $y=y_1-d_1(t)$,  and 
$$ 
B_1(t)=\left( 
\ba{cccc} 
0 & 0 & 0 & 0 \\ 
0 & 0 & 0 & 0 \\ 
0 & 0 & 0 & B_{v(t)}\!-\!B_{v(t_1)} \\ 
0 & 0 & \langle\na(\psi_{v(t)}\!-\!\psi_{v(t_1)}),\na\rho\rangle & 0 
\ea 
\right). 
$$ 
At the end of this section, we will derive appropriate bounds for the 
``remainder terms'' $B_1(t)Z_1(t)$ and $N_1(t)$ in (\re{redy1}). 
First, note that we have by Lemma \re{mod}, 
\be\la{BB1} 
|B_{v(t)}-B_{v(t_1)}|\le|\int_{t_1}^t\dot v(s)\cdot\na_vB_{v(s)}ds| 
\le C\int_{t}^{t_1}\Vert Z(s)\Vert_{-\beta}^2ds. 
\ee 
Similarly, 
\be\la{psipsi1} 
|\langle\na(\psi_{v(t)}-\psi_{v(t_1)}),\na\rho\rangle|\le 
 C\int_t^{t_1}\Vert Z(s)\Vert_{-\beta}^2ds. 
\ee 
Let us recall the following 
well-known inequality: 
for any $\al\in\R$ 
\be\la{pitre} 
(1+|y+x|)^{\alpha}\le(1+|y|)^{\alpha}(1+|x|)^{|\alpha|}, 
\,\,\,~~~~~~x,y\in\R^3. 
\ee 
\begin{lemma}\la{dest} 
For $(\Psi,\Pi,Q,P)\in{\cal E}_{\alpha}$ with any $\alpha\in\R$ 
the following estimate holds: 
\be\la{shiftest} 
\Vert(\Psi(y_1-d_1),\Pi(y_1-d_1),Q,P)\Vert_{\alpha}\le 
\Vert(\Psi,\Pi,Q,P)\Vert_{\alpha}(1+|d_1|)^{|\alpha|}~,\,\,\,~~~~~~d_1\in\R^3. 
\ee 
\end{lemma} 
{\bf Proof } Let us check the estimate only for one component, say, 
for $\Pi$. One has by (\re{pitre}) 
$$ 
\Vert\Pi(y_1-d_1,t)\Vert_{0,\alpha}^2=\int |\Pi(y_1-d_1,t)|^2 
(1+|y_1|)^{2\alpha}dy_1=\int |\Pi(y,t)|^2(1+|y+d_1|)^{2\alpha}dy\le 
$$ 
$$ 
\int|\Pi(y,t)|^2(1+|y|)^{2\alpha}(1+|d_1|)^{2|\alpha|}dy\le 
(1+|d_1|)^{2|\alpha|}\Vert\Pi\Vert^2_{0,\alpha}~, 
$$ 
and the lemma is proved.\hfill$\bo$ 
\begin{cor}\la{cor1} 
The following  bound holds 
\be\la{N1est} 
\Vert N_1(t)\Vert_{\beta}\le 
(1+|d_1(t)|)^{\beta}\Vert Z(t)\Vert^2_{-\beta}~,~~~~~~0\le t\le t_1. 
\ee 
\end{cor} 
Indeed, applying the previous lemma, we obtain from (\re{redN}) that 
$$ 
\Vert N_1(t)\Vert_{\beta}\le(1+|d_1(t)|)^{\beta} 
\Vert \ti N(t,Z(t))\Vert_{\beta}\le 
(1+|d_1(t)|)^{\beta}\Vert Z(t)\Vert^2_{-\beta}. 
$$ 
\begin{cor}\la{cor2} 
The following bound holds 
\be\la{B1Z1est} 
\Vert B_1(t)Z_1(t)\Vert_{\beta}\le C 
\Vert Z(t)\Vert_{-\beta}\int_t^{t_1} 
\Vert Z(\tau)\Vert^2_{-\beta} 
d\tau~,~~~~~~0\le t\le t_1. 
\ee 
\end{cor} 
For the proof we apply Lemma \re{dest} to (\re{BB1}) and (\re{psipsi1}) 
and use the fact that $B_1(t)Z_1(t)$ depends only on the finite-dimensional 
components of $Z_1(t)$. 
 
 
\setcounter{equation}{0}

\section{Integral Inequality}

Equation (\re{redy1}) can be represented in the integral form: 
\be\la{Z1duh} 
Z_1(t)=e^{A_1t}Z_1(0)+\int_0^te^{A_1(t-s)}[B_1Z_1(s)+N_1(s)]ds,\,\,\, 
0\le t\le t_1. 
\ee 
We apply the symplectic orthogonal 
projection $\bP_1:=\bP_{v(t_1)}$ to both sides, and get 
$$ 
\bP_1Z_1(t)=e^{A_1t}\bP_1Z_1(0)+\int_0^te^{A_1(t-s)}\bP_1[B_1Z_1(s)+N_1(s)]ds. 
$$ 
We have used here that $\bP_1$ commutes 
with 
the group $e^{A_1t}$ since the space $\cZ_1:=\bP_1\cE$ is invariant 
with respect to $e^{A_1t}$ by 
Remark \re{rZ}. 
Applying (\re{frozenest}) 
 we obtain  
\be\la{bPZ} 
\Vert \bP_1Z_1(t)\Vert_{-\beta}\le\fr{C}{(1+t)^{3/2}} 
\Vert \bP_1Z_1(0)\Vert_{\beta}+C\int_0^t\fr1{(1+|t-s|)^{3/2}}\Vert 
 \bP_1[B_1Z_1(s)+N_1(s)]\Vert_{\beta}ds. 
\ee 
The operator $\bP_1=\bI-\bPi_1$ is continuous in $\cE_\beta$ by (\re{Piv}). 
Hence,  from (\re{bPZ}) and 
(\re{N1est}), (\re{B1Z1est}), we obtain  
\beqn\la{duhest} 
\!\!\!\!\!\!\!\!\!\!\!\!\!\!\!\!\!\!&&\!\!\!\!\!\! 
\Vert \bP_1Z_1(t)\Vert_{-\beta}\le 
\fr{C(\ov d_1(0))}{(1+t)^{3/2}}\Vert Z(0)\Vert_{\beta} 
\nonumber\\ 
\!\!\!\!\!\!\!\!\!\!\!\!\!\!\!\!\!\!&&\!\!\!\!\!\! 
+C(\ov d_1(t))\int_0^t\fr1{(1+|t-s|)^{3/2}}\left[\Vert Z(s)\Vert_{-\beta} 
\int_s^{t_1}\Vert Z(\tau)\Vert^2_{-\beta}d\tau+ 
\Vert Z(s)\Vert^2_{-\beta}\right]ds,\,\,\,0\le t\le t_1, 
\eeqn 
where $\ov d_1(t):=\sup_{0\le s\le t} |d_1(s)| $. 
 
\begin{definition} Let $t_{*}'$ be the exit time 
\be\la{t*'} 
t_*'=\sup \{t\in[0,t_*): 
\ov d_1(s)\le 1,~~0\le s\le t\}. 
\ee 
\end{definition} 
Now (\re{duhest}) implies that for $t_1<t_*'$ 
\beqn\la{duhestri} 
\!\!\!\!\!\!\!\!\!\!\!\!\!\!\!\!\!\!&&\!\!\!\!\!\! 
\Vert \bP_1Z_1(t)\Vert_{-\beta}\le\fr{C} 
{(1+t)^{3/2}}\Vert Z(0)\Vert_{\beta} 
\nonumber\\ 
\!\!\!\!\!\!\!\!\!\!\!\!\!\!\!\!\!\!&&\!\!\!\!\!\! 
+C_1\int_0^t\fr1{(1+|t-s|)^{3/2}}\left[\Vert Z(s)\Vert_{-\beta} 
\int_s^{t_1}\Vert Z(\tau)\Vert^2_{-\beta}d\tau+ 
\Vert Z(s)\Vert^2_{-\beta}\right]ds,\,\,\,0\le t\le t_1. 
\eeqn 
 
 
\setcounter{equation}{0} 
 
\section{Symplectic Orthogonality}

Finally, we are going to change $\bP_1Z_1(t)$ by $Z(t)$ in the 
left hand side of 
(\re{duhestri}). 
We shall prove that this change is possible indeed by using 
    again the smalness condition  (\re{close}).  
For the justification we 
reduce the exit time  further. 
First, introduce the ``majorant'' 
\be\la{maj} 
m(t):= 
\sup_{s\in[0,t]}(1+s)^{3/2}\Vert Z(s)\Vert_{-\beta}~,~~~~~~~~~t\in [0,t_*). 
\ee 
Denote by $\ve$ a fixed positive number (which  
will be specified below). 
\begin{definition} Let $t_{*}''$ be the exit time 
\be\la{t*''} 
t_*''=\sup \{t\in[0,t_*'): 
m(s)\le \ve,~~0\le s\le t\}. 
\ee 
\end{definition}

The following important bound  (\re{Z1P1est}) 
enables us to change the norm of $\bP_1Z_1(t)$ on the left hand side 
 of 
(\re{duhestri}) by the norm of $Z(t)$.

\begin{lemma}\la{Z1P1Z1} 
For sufficiently small $\ve>0$, we have  
\be\la{Z1P1est} 
\Vert Z(t)\Vert_{-\beta}\le C\Vert \bP_1Z_1(t)\Vert_{-\beta}, 
~~~~~~~~0\le t \le t_1, 
\ee 
for any $t_1<t_*''$, where $C$ depends only on $\rho$ and $\ov v$. 
\end{lemma} 
{\bf Proof } 
The proof is based on 
the symplectic orthogonality 
(\re{PZ}), i.e. 
\be\la{PZ1} 
\bPi_{v(t)}Z(t)=0,~~~~t\in[0,t_1], 
\ee 
and on the fact 
that 
 all the spaces $\cZ(t):=\bP_{v(t)}\cE$ are almost parallel 
for all $t$ (see Fig. 2). 
\begin{figure}[htbp] 
\begin{center} 
\includegraphics[width=0.8\columnwidth]{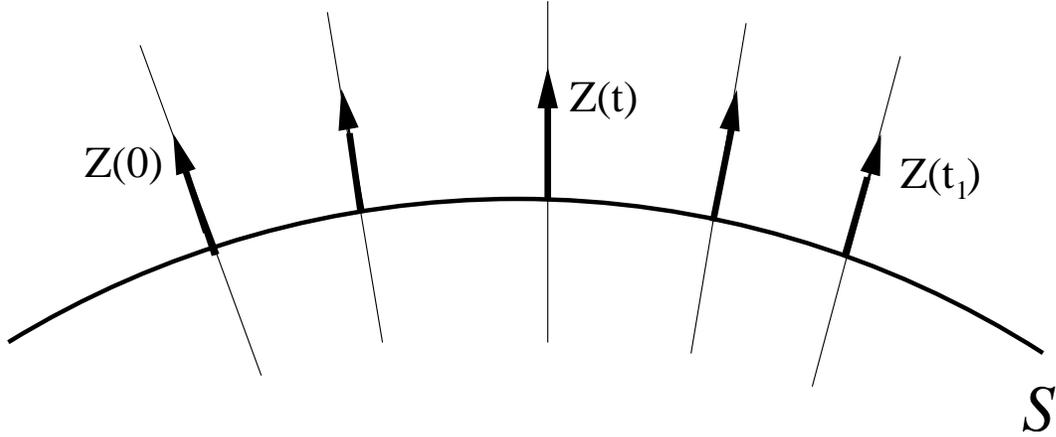} 
\caption{Symplectic orthogonality.} 
\label{fig-2} 
\end{center} 
\end{figure} 
Namely, we first note that 
$\Vert Z(t)\Vert_{-\beta}\le C\Vert Z_1(t)\Vert_{-\beta}$ 
by Lemma \re{dest}, since $|d_1(t)|\le 1$ for 
 $t\le t_1<t_*''<t_*'$. Therefore, it suffices 
 to prove that 
\be\la{Z1P1ests} 
\Vert Z_1(t)\Vert_{-\beta}\le 2\Vert \bP_1Z_1(t)\Vert_{-\beta}, 
~~~~~~~~0\le t\le t_1. 
\ee 
This estimate will follow from 
\be\la{Z1P1estf} 
\Vert\bPi_{v(t_1)}Z_1(t)\Vert_{-\beta}\le\fr12\Vert Z_1(t)\Vert_{-\beta}, 
\,\,\,0\le t\le t_1\,, 
\ee 
since $\bP_1Z_1(t)=Z_1(t)-\bPi_{v(t_1)}Z_1(t)$. 
To prove (\re{Z1P1estf}), we write (\re{PZ1}) as 
\be\la{PZ0r} 
\bPi_{v(t),1}Z_1(t)=0,~~~~t\in[0,t_1], 
\ee 
where $\bPi_{v(t),1}Z_1(t)$ is $\bPi_{v(t)}Z(t)$ expressed in 
terms of the variable $y_1=y+d_1(t)$. 
Hence, (\re{Z1P1estf}) follows from  (\re{PZ0r}) if 
the difference 
 $\bPi_{v(t_1)}-\bPi_{v(t),1}$ is small uniformly in $t$, 
i.e. 
\be\la{difs} 
\Vert\bPi_{v(t_1)}-\bPi_{v(t),1}\Vert<1/2,~~~~~~~0\le t\le t_1. 
\ee 
 
It remains to justify (\ref{difs}) for any  
sufficiently small $\varepsilon >0$. 
We will need the formula 
(\ref{Piv}) and the following relation which follows from 
(\ref{Piv}): 
\begin{equation} 
\mbox{\boldmath $\Pi$}_{v(t),1}Z_{1}(t)=\sum \mbox{\boldmath $\Pi$}%
_{jl}(v(t))\tau _{j,1}(v(t))\Omega (\tau _{l,1}(v(t)),Z_{1}(t)),  \label{P1} 
\end{equation} 
where $\tau _{j,1}(v(t))$ are the vectors $\tau _{j}(v(t))$  
expressed via the 
variables $y_{1}$. In detail (cf. (\ref{inb})), 
\begin{equation} 
\left. 
\begin{array}{rcl} 
\tau _{j,1}(v) & := & (-\partial _{j}\psi _{v}(y_{1}-d_{1}(t)),-\partial 
_{j}\pi _{v}(y_{1}-d_{1}(t)),e_{j},0), \\ 
\tau _{j+3,1}(v) & := & (\partial _{v_{j}}\psi _{v}(y_{1}-d_{1}(t)),\partial 
_{v_{j}}\pi _{v}(y_{1}-d_{1}(t)),0,\partial _{v_{j}}p_{v}), 
\end{array} 
\right| ~~~j=1,2,3,  \label{inb*} 
\end{equation} 
where $v=v(t)$. Since $|d_{1}(t)|\leq 1$, and the functions 
$\nabla \tau _{j}$  
are smooth and rapidly 
decaying at infinity, Lemma \re{dest} 
implies that 
\begin{equation} 
\Vert \tau _{j,1}(v(t))-\tau _{j}(v(t))\Vert_\beta\leq 
C|d_{1}(t)|^\beta,~~~~0\leq t\leq t_{1}  \label{011} 
\end{equation} 
for all $j=1,2,\dots,6$. 
Furthermore, 
\[ 
\tau _{j}(v(t))-\tau _{j}(v(t_{1}))=\int_{t}^{t_{1}}\dot{v}(s)\cdot \nabla 
_{v}\tau _{j}(v(s))ds, 
\] 
and therefore 
\begin{equation} 
\Vert \tau _{j}(v(t))-\tau _{j}(v(t_{1}))\Vert_\beta\leq 
C\int_{t}^{t_{1}}|\dot{v}(s)|ds,~~~~0\leq t\leq t_{1}.  \label{012} 
\end{equation} 
Similarly, 
\begin{equation} 
|\mbox{\boldmath $\Pi$}_{jl}(v(t))-\mbox{\boldmath 
$\Pi$}_{jl}(v(t_{1}))|= |\int_{t}^{t_{1}}\dot{v}(s)\cdot \nabla 
_{v}\mbox{\boldmath $\Pi$}_{jl}(v(s))ds|\leq 
C\int_{t}^{t_{1}}|\dot{v}(s)|ds,~~~~0\leq t\leq t_{1}, \label{013} 
\end{equation} 
since $|\nabla_{v}\mbox{\boldmath $\Pi$}_{jl}(v(s))|$ is uniformly bounded by (\re{sit}). 
Hence, the bounds (\ref{difs}) will follow from (\ref{Piv}), 
(\ref{P1}) and (\ref{011})-(\ref{013}) if we  
shall prove 
 that 
$|d_{1}(t)|$ and the integral on the right hand side of 
(\ref{012}) can be made as small as desired by choosing 
a sufficiently small $\varepsilon >0$. 
 
To estimate $d_{1}(t)$, note that 
\be\la{wen} 
w(s)-v(t_{1})=w(s)-v(s)+v(s)-v(t_{1})= 
\dot c(s)+\int_s^{t_1}\dot v(\tau)d\tau 
\ee 
by (\re{vw}). 
Hence, the definitions   
(\ref {dd1}), (\re{maj}), and   
Lemma \ref{mod} imply that 
\begin{eqnarray} 
\!\!\!\!\!\! \!\!\! \!\!\! 
 |d_{1}(t)|\!\!\! &=&\!\!\!|\int_{t_{1}}^{t}(w(s)-v(t_{1}))ds|\leq 
\int_{t}^{t_{1}}\left( |\dot{c}(s)|+\int_{s}^{t_{1}}|\dot{v}(\tau )|d\tau 
\right)ds   \nonumber  \label{d1est} \\ 
&&  \nonumber \\ 
\!\!\! &\leq &\!\!\!Cm^{2}(t_{1})\int_{t}^{t_{1}}\left( \frac{1}{(1+s)^{3}}%
+\int_{s}^{t_{1}}\frac{d\tau }{(1+\tau )^{3}}\right) ds\leq 
Cm^{2}(t_{1})\leq C\varepsilon ^{2},~~~~0\leq t\le t_{1} 
\end{eqnarray} 
since $t_{1}<t_{\ast }^{\prime \prime }$. Similarly, 
\begin{equation} 
\int_{t}^{t_{1}}|\dot{v}(s)|ds\leq Cm^{2}(t_{1})\int_{t}^{t_{1}}\frac{ds}{%
(1+s)^{3}}\leq C\varepsilon ^{2},~~~~0\leq t\le t_{1}. 
\label{tvjest} 
\end{equation} 
The proof is completed. 
\bo

 
\setcounter{equation}{0} 
 
\section{Decay of Transversal Component} 
Here we prove Proposition \re{pdec}. 
\\ 
{\it Step i)} We fix $\ve>0$ and $t''_*=t''_*(\ve)$ 
for which Lemma \re{Z1P1Z1} holds. 
Then a bound of type (\re{duhestri}) 
holds with 
$\Vert \bP_1Z_1(t)\Vert_{-\beta}$  
replaced by 
 $\Vert Z(t)\Vert_{-\beta}$ 
on the left hand side: 
\beqn\la{duhestrih} 
\!\!\!\!\!\!\!\!\!\!\!\!\!\!\!\!\!\!&&\!\!\!\!\!\! 
\Vert Z(t)\Vert_{-\beta}\le\fr{C} 
{(1+t)^{3/2}}\Vert Z(0)\Vert_{\beta} 
\nonumber\\ 
\!\!\!\!\!\!\!\!\!\!\!\!\!\!\!\!\!\!&&\!\!\!\!\!\! 
+C\int_0^t\fr1{(1+|t-s|)^{3/2}}\left[\Vert Z(s)\Vert_{-\beta} 
\int_s^{t_1}\Vert Z(\tau)\Vert^2_{-\beta}d\tau+ 
\Vert Z(s)\Vert^2_{-\beta}\right]ds,\,\,\,0\le t\le t_1 
\eeqn 
for $t_1<t_*'$. 
This implies an integral inequality 
 for the majorant 
$$ 
m(t):=\sup_{s\in[0,t]}(1+s)^{3/2}\Vert Z(s)\Vert_{-\beta}~. 
$$ 
Namely, multiplying  (\re{duhestrih}) 
by $(1+t)^{3/2}$ 
and taking the supremum in $t\in[0,t_1]$, 
we get 
$$ 
m(t_1) 
\le C\Vert Z(0)\Vert_{\beta}+ 
C\sup_{t\in[0,t_1]}\ds 
\int_0^t\fr{(1+t)^{3/2}}{(1+|t-s|)^{3/2}}\left[\fr{m(s)}{(1+s)^{3/2}} 
\int_s^{t_1}\fr{m^2(\tau)d\tau}{(1+\tau)^{3}}+\fr{m^2(s)} 
{(1+s)^{3}}\right]ds 
$$ 
for $t_1< t_*''$. Taking into account that $m(t)$ is a monotone 
increasing function, we get \be\la{mest} m(t_1)\le C\Vert 
Z(0)\Vert_{\beta}+C[m^3(t_1)+m^2(t_1)]I(t_1), ~~~~~~~~~~~~~t_1< 
t_*''. \ee where 
$$ 
I(t_1)= 
\sup_{t\in[0,t_1]} 
\int_0^{t}\fr{(1+t)^{3/2}}{(1+|t-s|)^{3/2}}\left[\fr1{(1+s)^{3/2}} 
\int_s^{t_1}\fr{d\tau}{(1+\tau)^{3}}+\fr1{(1+s)^3}\right]ds 
\le \ov I<\infty,~~~~t_1\ge0. 
$$ 
Therefore, (\re{mest}) becomes 
\be\la{m1est} 
m(t_1)\le C\Vert Z(0)\Vert_{\beta}+C\ov I[m^3(t_1)+m^2(t_1)],~~~~ t_1<t_*''. 
\ee 
This inequality implies that $m(t_1)$ 
is bounded for $t_1<t_*''$, and moreover, 
\be\la{m2est} 
m(t_1)\le C_1\Vert Z(0)\Vert_{\beta},~~~~~~~~~t_1<t_*''\,, 
\ee 
since $m(0)=\Vert Z(0)\Vert_{\beta}$ is sufficiently small by 
 (\re{closeZ}). 
\\ 
{\it Step ii)} The constant $C_1$ in the estimate 
(\re{m2est}) does not depend on 
$t_*$, $t_*'$ and $t_*''$ by Lemma \re{Z1P1Z1}. 
We choose a small $d_0$ in (\re{close}) such that 
$\Vert Z(0)\Vert_{\beta}<\ve/(2C_1)$.  
This is possible by (\re{closeZ}). 
Then the estimate (\re{m2est}) implies that $t''_*=t'_*$~,  
and therefore 
 (\re{m2est}) holds for all $t_1<t'_*$. 
Then the bound (\re{d1est}) holds for all $t<t_*'$. Choose  
a small $\ve$ such 
that the right hand side in (\re{d1est})  
does not exceed one. Then $t'_*=t_*$. 
Therefore, 
(\re{m2est}) holds for any $t_1<t''_*=t_*$, hence 
(\re{Zt}) also holds if 
$\Vert Z(0)\Vert_{\beta}$ is sufficiently small. 
Finally, this implies that 
$t_*=\infty$. Hence we also have $t''_*=t'_*=\infty$, and 
(\re{m2est}) holds for any $t_1>0$ if $d_0$ is sufficiently  
small. 
 \bo 
 
\section{Soliton Asymptotics} 
\setcounter{equation}{0}

Here we prove our main 
Theorem \re{main} under the assumption 
that 
the decay (\re{Zdec}) holds. Let us first prove the asymptotics 
(\re{qq}) for the vector components, and then 
the asymptotics (\re{S}) for the fields. 
\\ 
{\bf Asymptotics for the vector components} 
It follows from (\re{addeq}) that $\dot q=\dot b+\dot Q$, and from (\re{reduced}), 
(\re{redN}), and (\re{AA}) that $\dot Q=B_{v(t)}P+{\cal O} 
(\Vert Z\Vert^2_{-\beta})$. Thus, 
\be\la{dq} 
\dot q=\dot b+\dot Q=v(t)+\dot c(t)+B_{v(t)}P(t)+{\cal O} 
(\Vert Z\Vert^2_{-\beta}). 
\ee 
Equation 
(\re{parameq}) together with estimates (\re{NZ}) and (\re{Zdec}) 
imply that 
\be\la{bv} 
|\dot c(t)|+|\dot v(t)|\le \ds\fr {C_1(\rho,\ov v,d_0)}{(1+t)^{3}}, 
~~~~~~t\ge0. 
\ee 
Therefore, $c(t)=c_+ +\cO(t^{-2})$ and $v(t)=v_+ +\cO(t^{-2})$, $t\to\infty$. 
Since $|P|\le\Vert Z\Vert_{-\beta}$, the estimate (\re{Zdec}) 
together with relations (\re{bv}) and (\re{dq}) imply that 
\be\la{qbQ} 
\dot q(t)=v_++\cO(t^{-3/2}). 
\ee 
Similarly, 
\be\la{bt} 
b(t)=c(t)+\ds\int_0^tv(s)ds=v_+t+a_++\cO(t^{-1}), 
\ee 
and hence 
the second part of (\re{qq}) follows: 
\be\la{qbQ2} 
q(t)=b(t)+Q(t)=v_+t+a_++\cO(t^{-1}), 
\ee 
since $Q(t)=\cO(t^{-3/2})$ 
by  (\re{Zdec}). 
\\ 
{\bf Asymptotics for the fields} We apply the approach 
developed in \ci{IKSs,KKS}. 
For the field part of the solution, $F(t)=(\psi(x,t),\pi(x,t))$, 
in the original variable $x$, 
let us define the accompanying soliton field as 
$F_{\rm v(t)}(t)=(\psi_{\rm v(t)}(x-q(t)),\pi_{\rm v(t)}(x-q(t)))$, 
where we now set ${\rm v}(t)=\dot q(t)$, cf. (\re{dq}). 
Then for the difference $Z(t)=F(t)-F_{\rm v(t)}(t)$ we obtain 
easily, from the first two equations of the system 
(\re{system0}), 
 the inhomogeneous Klein-Gordon equation \ci[(2.5)]{KKS}, 
$$ 
\dot Z(t)=A_0Z(t)-\dot{\rm v}\cdot\na_{\rm v}F_{{\rm v}(t)}(t), 
\,\,\,\,\,\,A_0(\psi,\pi)=(\pi,(\De-m^2)\psi). 
$$
Then 
\be\la{eqacc} 
Z(t)=W_0(t)Z(0)-\int_0^tW_0(t-s)[\dot{\rm v}(s)\cdot\na_{\rm v} 
F_{{\rm v}(s)}(s)]ds, 
\ee 
where $W_0(t)$ is the dynamical group of free Klein-Gordon equation. 
To obtain the asymptotics (\re{S}) it suffices to prove that 
$Z(t)=W_0(t)\bPsi_++r_+(t)$ for some $\bPsi_+\in{\cal F}$ and 
that  
$\Vert r_+(t)\Vert_{{\cal F}}={\cal O}(t^{-1/2})$. 
This is equivalent to the asymptotics 
\be\la{Sme} 
W_0(-t)Z(t)=\bPsi_++r_+'(t),~~~~~ 
\Vert r_+'(t)\Vert_\cF=\cO(t^{-1/2}), 
\ee 
since 
$W_0(t)$ is a unitary group on the Sobolev space $\cF$ by the energy 
conservation for the free Klein-Gordon equation. Finally, 
the asymptotics 
 (\re{Sme}) hold since 
 (\re{eqacc}) implies that 
\be\la{duhs} 
W_0(-t)Z(t)= 
Z(0)-\int_0^t W_0(-s)R(s)ds,\,\,\,\,\,R(s) 
=\dot{\rm v}(s)\cdot\na_{\rm v}F_{{\rm v}(s)}(s), 
\ee 
where 
the integral  on the right hand side of (\re{duhs}) 
converges in the Hilbert space $\cF$ 
with the rate $\cO(t^{-1/2})$. The latter holds since 
$\Vert  W_0(-s)R(s)\Vert_\cF =\cO(s^{-3/2})$ by the unitarity 
of $W_0(-s)$ and the decay rate $\Vert R(s)\Vert_\cF =\cO(s^{-3/2})$. 
Let us prove this rate of decay. It suffices to prove that 
$|\dot {\rm v}(s)|=\cO(s^{-3/2})$, or equivalently 
$|\dot p(s)|=\cO(s^{-3/2})$. Substitute (\re{add}) to the last 
equation of (\re{system0}) and obtain 
\beqn 
\dot p(t)&=&\int\left[\psi_{v(t)}(x-b(t))+\Psi(x-b(t),t)\right] 
\na\rho(x-b(t)-Q(t))dx 
\bigskip\\ 
&=&\int\psi_{v(t)}(y)\na\rho(y)dy+ 
\int\psi_{v(t)}(y)\left[\na\rho(y-Q(t))-\na\rho(y)\right]dy 
+\int\Psi(y,t)\na\rho(y-Q(t))dy.\nonumber 
\eeqn 
The first integral on the right hand side is zero by the stationary 
equations (\re{stfch}). The second integral is ${\cal O}(t^{-3/2})$, 
which follows  
from  the conditions (\re{ro}) on $\rho$ 
and the asymptotics 
 $Q(t)={\cal O}(t^{-3/2})$. 
Finally, the third integral is ${\cal O}(t^{-3/2})$ 
by estimate (\re{Zdec}). This completes the proof.\bo 
 
 
\setcounter{equation}{0} 
 
\section{Decay for the Linearized Dynamics} 
In remaining section, we 
prove Proposition \re{lindecay} 
to complete the proof of the main result (Theorem \re{main}). 
Here we discuss the general strategy of proving the Proposition. 
We apply the Fourier-Laplace transform 
\be\la{FL} 
\ti X(\lam)=\int_0^\infty e^{-\lam t}X(t)dt,~~~~~~~\Re\lam>0 
\ee 
to (\re{Avv}).  According to Proposition \re{lindecay}, we can 
expect that 
the solution $X(t)$ is bounded in the norm $\Vert\cdot\Vert_{-\beta}$. 
Then the integral (\re{FL}) converges and is analytic for $\Re\lam>0$, 
and 
\be\la{PW} 
\Vert\ti X(\lam)\Vert_{-\beta}\le \ds\fr{C}{\Re\lam},~~~~~~~\Re\lam>0. 
\ee 
Let us derive an equation for $\ti X(\lam)$ which is equivalent to the 
Cauchy problem for (\re{Avv}) with the initial condition $X(0)=X_0\in\cE_{-\beta}$. 
We shall write $A$ and $v$ instead of $A_1$ and $v_1$ 
in all remaining part of the paper. 
Applying the Fourier-Laplace transform to (\re{Avv}), we get 
that 
\be\la{FLA} 
\lam\ti X(\lam)=A\ti X(\lam)+X_0,~~~~~~~~\Re\lam>0. 
\ee 
Let us stress that 
 (\re{FLA}) 
 is equivalent to 
the Cauchy problem for the functions $X(t)\in 
C_b([0,\infty);\cE_{-\beta})$. 
Hence the solution $X(t)$ is given by 
\be\la{FLAs} 
\ti X(\lam)=-(A-\lam)^{-1}X_0,~~~~~~~~\Re\lam>0 
\ee 
if the resolvent $R(\lam)=(A-\lam)^{-1}$ exists for $\Re\lam>0$. 
 
Let us comment on our following strategy in proving  the decay 
(\re{Zdec}). 
We shall first 
construct the resolvent  $R(\lam)$ for $\Re\lam>0$ 
and prove that  
this resolvent 
is a continuous operator on ${\cal E}_{-\beta}$. 
Then $\ti X(\lam)\in\cE_{-\beta}$ and 
is an analytic function for $\Re\lam>0$. 
After this we must to justify that there exists a (unique) function 
$X(t)\in C([0,\infty);\cE_{-\beta})$ satisfying (\re{FL}).

The analyticity of $\ti X(\lam)$ and the 
Paley-Wiener arguments (see \ci{EKS}) should provide 
the existence of a $\cE_{-\beta}$ - valued distribution $X(t)$, $t\in\R$, 
with a support in $[0,\infty)$. Formally, 
\be\la{FLr} 
X(t)=\fr1{2\pi}\int_\R e^{i\om t}\ti X(i\om+0)d\om, ~~~~~~~~t\in\R. 
\ee 
However, to establish the continuity of $X(t)$ for $t\ge 0$, 
we need an additional bound for $\ti X(i\om+0)$ for large  
values of $|\om|$. 
Finally, for the time decay of $X(t)$, we need an additional 
information on the 
smoothness and decay of $\ti X(i\om+0)$. 
More precisely, 
we must prove that the function $\ti X(i\om+0)$  
has the following properties: 
\\ 
 i)  it is smooth outside 
$\om=0$ and $\om=\pm\mu$, where $\mu=\mu(v)>0$, 
\\ 
ii) it decays in a certain sense as $|\om|\to\infty$, 
\\ 
iii) 
it admits the Puiseux expansion at $\om=\pm\mu$, 
\\ 
iv) it is analytic at $\om=0$ if  
$X_0\in\cZ_v:=\bP_v\cE$ and $X_0\in\cE_\beta$. 
\\ 
Then the decay (\re{Zdec}) would follow from the Fourier-Laplace 
representation (\re{FLr}).

We shall check with detail 
properties of the type 
i)-iv) only for the last two  components 
$\ti Q(\lam)$ and $\ti P(\lam)$ of the vector $\ti X(\lam) 
=(\ti\Psi(\lam),\ti\Pi(\lam),\ti Q(\lam),\ti P(\lam))$. 
The properties 
provide the decay (\re{Zdec}) for the vector components 
$Q(t)$ and $P(t)$ of the solution $X(t)$. 
 
However, we will not prove the properties 
 of the type i)-iv) for the field components 
$\Psi(x,\lam)$ and $\Pi(x,\lam)$. 
The decay (\re{Zdec}) for the field components 
is deduced 
in Section 18 
directly from the time-dependent 
field equations of the system (\re{Avv}), 
using the decay of the component $Q(t)$ and a 
version of  strong Huygens principle for the Klein-Gordon equation.

 
 
\setcounter{equation}{0} 
 
\section{Constructing the Resolvent} 
 
To justify the representation (\re{FLAs}), 
we 
construct the resolvent as a bounded operator in 
${\cal E}_{-\beta}$ 
 for $\Re\lam>0$. 
We shall write $(\Psi(y),\Pi(y), Q, P)$ 
instead of $(\ti\Psi(y,\lam),\ti\Pi(y,\lam),\ti Q(\lam),\ti P(\lam))$ 
to simplify the notations. Then (\re{FLA}) reads 
$$ 
(A-\lam)\left( 
\ba{c} 
\Psi \\ \Pi \\ Q \\ P 
\ea 
\right)=-\left( 
\ba{c} 
\Psi_0 \\ \Pi_0 \\ Q_0 \\ P_0 
\ea 
\right),\,\,\,{\rm where}\,\,\,A\left( 
\ba{c} 
\Psi \\ \Pi \\ Q \\ P 
\ea 
\right)=\left( 
\ba{r} 
\Pi+v\cdot\na\Psi \\ 
\De\Psi-m^2\Psi+v\cdot\na\Pi+Q\cdot\na\rho \\ 
B_vP \\ 
-\langle\na\Psi,\rho\rangle+\langle\na\psi_{v},Q\cdot\na\rho\rangle 
\ea 
\right). 
$$ 
This gives  the system of equations 
\be\la{eq1} 
\left.\ba{r} 
\Pi(y)+v\cdot\na\Psi(y)-\lam\Psi(y)=-\Psi_0 (y) 
\\ 
\\ 
\De\Psi(y)-m^2\Psi(y)+v\cdot\na\Pi(y)+Q\cdot\na\rho(y)-\lam\Pi(y)=-\Pi_0(y) 
\\ 
\\ 
B_vP-\lam Q=-Q_0 
\\ 
\\ 
-\langle\na\Psi(y),\rho(y)\rangle+\langle\na\psi_{v}(y),Q\cdot\na\rho(y) 
\rangle-\lam P=-P_0 
\ea\right|~~~~~~~~~y\in\R^3. 
\ee 
{\it Step i)} Let us study the first two equations. In the  
Fourier space 
 they become 
\be\la{F1} 
\left. 
\ba{rcl} 
\hat\Pi(k)-ivk\hat\Psi(k)-\lam\hat\Psi(k)&=&-\hat\Psi_0(k) 
\\\\ 
(-k^2-m^2)\hat\Psi(k)-(ivk+\lam)\hat\Pi(k)&=&-\hat\Pi_0(k)+iQk\hat\rho(k) 
\ea\right|~~~~~~~~~k\in\R^3\5. 
\ee 
Let us invert the matrix of the system and obtain 
$$ 
\left( 
\ba{cc} 
-(ivk+\lam) & 1 \\ 
-(k^2+m^2) & -(ivk+\lam) 
\ea 
\right)^{-1}=[(ivk+\lam)^2+k^2+m^2]^{-1}\left( 
\ba{cc} 
-(ivk+\lam) & -1 \\ 
k^2+m^2 & -(ivk+\lam) 
\ea 
\right). 
$$ 
Taking the inverse Fourier transform, 
we obtain the corresponding fundamental solution 
\be\la{Green} 
G_{\lam}(y)=\left( 
\ba{cc} 
v\cdot\na-\lam & -1 \\ 
-\De+m^2 & v\cdot\na-\lam 
\ea 
\right)g_{\lam}(y), 
\ee 
where $g_{\lam}(y)$ is the 
unique tempered 
fundamental solution of the determinant 
\be\la{fso} 
D=D(\lam)=-\De+m^2+(-v\cdot\na+\lam)^2. 
\ee 
From now on we use the system of coordinates in $x$-space in which 
$v=(|v|,0,0)$, hence $vk=|v|k_1$, and 
\be\la{dete} 
g_\lam(y)=F^{-1}_{k\to y}\ds\fr{1}{k^2+m^2+(ivk+\lam)^2}= 
F^{-1}_{k\to y}\ds\fr{1}{k^2+m^2+ 
(i|v|k_1+\lam)^2},~~~y\in\R^3. 
\ee 
Note that the denominator does not vanish for $\Re\lam>0$. 
This implies 
\bl\la{cres} 
The operator $G_{\lam}$ with the integral kernel $G_{\lam}(y-y')$ 
is continuous as an operator from 
$H^1(\R^3)\oplus L^2(R^3)$ to $H^2(\R^3)\oplus H^1(R^3)$ 
 for $\Re\lam>0$. 
\el 
Thus, formulas (\re{F1}) and 
 (\re{Green}) imply the convolution representation 
\be\la{Psi} 
\left.\ba{rcl} 
\Psi&=&-(v\cdot\na-\lam)g_{\lam}*\Psi_0+g_{\lam}*\Pi_0+(g_{\lam}*\na\rho)\cdot Q 
\\\\ 
\Pi&=&-(-\De+m^2)g_{\lam}*\Psi_0- 
(v\cdot\na-\lam)g_{\lam}*\Pi_0-(v\cdot\na-\lam)(g_{\lam}*\na\rho)\cdot Q 
\ea\right| 
\ee 
\noindent{\it Step ii)} 
Let us compute $g_{\lam}(y)$ explicitly. 
First consider the case 
$v=0$. The fundamental solution of the operator $-\De+m^2+\lam^2$ is 
\be\la{fus} 
g_{\lam}(y)=\fr{e^{-\kappa|y|}}{4\pi|y|}, 
\ee 
where 
\be\la{kap} 
\kappa^2=m^2+\lam^2, ~~~~~~~~~~~ 
\Re\kappa>0 ~~\mbox{ for }~~ \Re\lam>0. 
\ee 
Thus, in the case $v=0$ we have 
$$ 
G_{\lam}(y-y')=\left( 
\ba{cc} 
-\lam & -1 \\ 
-\De+m^2 & -\lam 
\ea 
\right)\fr{e^{-\sqrt{\lam^2+m^2}|y-y'|}}{4\pi|y-y'|}. 
$$ 
For general $v=(|v|,0,0)$ 
with $|v|<1$ 
the denominator in (\re{dete}), which is the Fourier symbol of $D$, 
reads 
\beqn 
\hat D(k)&=& 
k^2+m^2+(i|v|k_1+\lam)^2 
\nonumber\\ 
 \nonumber\\ 
&=&(1-v^2)k_1^2+k_2^2+k_3^2+ 
2i|v|k_1\lam+\lam^2+m^2 
\nonumber\\ 
 \nonumber\\ 
&=& 
(1-v^2)(k_1+\fr{i|v|\lam}{1-v^2})^2+k_2^2+k_3^2+ \kappa^2, 
\la{Dhat} 
\eeqn 
where 
\be\la{tikappa} 
\kappa^2=\fr{v^2\lam^2}{1-v^2}+\lam^2+m^2=\fr{\lam^2}{1-v^2}+m^2. 
\ee 
Therefore, setting $\ga:=1/\sqrt{1-v^2}$, we have 
\be\la{kappa} \kappa=\ga\sqrt{\lam^2+\mu^2}, ~~~~~~\mu:=m/\ga. 
\ee 
Return to $x$-space: 
\be\la{13.10'} 
D=-\fr{1}{\ga^2}(\na_1+\ga\ka_1)^2- \na_2^2-\na_3^2+\kappa^2, 
~~~~~~~~ \ka_1:=\ga|v|\lam 
\ee 
Define 
$\ti y_1:=\ga y_1$ 
and $\ti\na_1:=\pa/\pa \ti y_1$. 
Then 
\be\la{Dy} 
D=-(\ti\na_1+\ka_1)^2-\na_2^2-\na_3^2+\kappa^2. 
\ee 
Thus, its fundamental solution is 
\be\la{glam} 
g_{\lam}(y)=\fr{e^{-\kappa|\ti y|-\ka_1\ti y_1}}{4\pi|\ti y|}, 
\,\,\,\ti y:=(\ga y_1,y_2,y_3), 
\ee 
where we choose $\Re\kappa>0$ for $\Re\lam>0$. 
Let us note that 
\be\la{kak} 
 0<\Re\ka_1<\Re\kappa,~~~~~~\Re\lam>0. 
\ee 
This inequality follows from 
the fact that 
the fundamental solution decays exponentially 
by the Paley-Wiener arguments since 
the quadratic form 
(\re{Dhat}) does not vanish in a complex neighborhood 
of the real space $\R^3$ for $\Re\lam>0$. 
Let us state the result which we have got above. 
\bl \la{cac} 
i) The 
operator $D=D(\lam)$ is invertible in $L^2(\R^3)$ for $\Re\lam>0$ and 
its fundamental solution (\re{glam}) decays exponentially. 
\\ 
ii) Formulas (\re{glam}) and  (\re{kappa}) imply that, 
for every fixed $y$, 
the Green function 
$g_\lam(y)$ admits an analytic continuation (in  
the variable $\lam$) 
to the 
Riemann surface of the 
algebraic function $\sqrt{\lam^2+\mu^2}$ with the 
branching points ~$\la=\pm i\mu$. 
\el 
\noindent{\it Step iii)} 
Let us now proceed with the last two equations (\re{eq1}), 
\be\la{lte} 
-\lam Q+B_{v}P=-Q_0,\,\,\,\langle\na\psi_{v}, 
Q\cdot\na\rho\rangle-\langle\na\Psi,\rho\rangle-\lam P=-P_0. 
\ee 
Let us eliminate the field $\Psi$ by the first equation 
(\re{Psi}). Namely, rewrite the equation in the form 
$\Psi(x)=\Psi_1(Q)+\Psi_2(\Psi_0,\Pi_0)$, 
where 
\be\la{Psi12} 
\Psi_1(Q) 
=Q\cdot(g_{\lam}*\na\rho), 
~~~~~~~~~ 
\Psi_2(\Psi_0,\Pi_0)=-(v\cdot\na-\lam)g_{\lam}*\Psi_0+g_{\lam}*\Pi_0. 
\ee 
Then we have 
$$ 
\langle\na\Psi,\rho\rangle=\langle\na\Psi_1(Q),\rho\rangle+ 
\langle\na\Psi_2(\Psi_0,\Pi_0),\rho\rangle, 
$$  
and the last equation in (\re{lte}) 
 becomes 
$$ 
\langle\na\psi_{v},Q\cdot\na\rho\rangle-\langle\na\Psi_1(Q), 
\rho\rangle-\lam P=-P_0+\langle\na\Psi_2(\Psi_0,\Pi_0),\rho\rangle=:-P_0'. 
$$ 
Let us first  
compute the term $\langle\na\psi_{v},Q\cdot\na\rho\rangle= 
\sum_{j} 
\langle\na\psi_{v},Q_j\pa_j\rho\rangle= 
\sum_{j} \langle\na\psi_{v},\pa_j\rho\rangle Q_j$. 
Applying the Fourier transform $F_{y\to k}$,  
the Parseval identity, and (\re{hpsiv}) we see that 
\beqn 
\!\!\!\!\!\!\!\! 
\langle\pa_i\psi_{v},\pa_j\rho\rangle&=& 
\langle -ik_i\hat\psi_{v}(k),-ik_j\hat\rho(k)\rangle= 
\langle k_i\hat\psi_{v}(k),k_j\hat\rho(k)\rangle= 
\nonumber\\ 
\nonumber\\ 
&&-\langle\fr{k_i\hat\rho(k)}{k^2+m^2-(|v|k_1)^2},k_j\hat\rho(k)\rangle= 
-\int\fr{k_ik_j|\hat\rho(k)|^2dk}{\ds k^2+m^2-(|v|k_1)^2}=:-K_{ij}. 
\la{Lij} 
\eeqn 
As the result, $\langle\na\psi_{v},Q\cdot\na\rho\rangle=-KQ$, where $K$ 
is the $3\times3$ matrix with the matrix elements $K_{ij}$. 
The matrix $K$ is diagonal and positive definite since $\hat\rho(k)$ 
 is spherically 
symmetric and not identically zero by (\re{W}).

Let us now  
compute the term $-\langle\na\Psi_1,\rho\rangle= 
\langle\Psi_1,\na\rho\rangle$. 
We have  
$$ 
\langle\Psi_1,\pa_i\rho\rangle= 
\langle\sum\limits_j(g_{\lam}*\pa_j\rho)Q_j,\pa_i\rho\rangle= 
\sum\limits_j\langle g_{\lam}*\pa_j\rho,\pa_i\rho\rangle Q_j= 
\sum\limits_j H_{ij}(\lam)Q_j 
$$ 
 since  
$\Psi_1=Q\cdot(g_{\lam}*\na\rho)$, 
and by the Parseval identity again, 
we have 
\beqn\la{Cij} 
H_{ij}(\lam):&=&\langle g_{\lam}*\pa_j\rho, 
\pa_i\rho\rangle= 
\langle i\hat g_{\lam}(k)k_j\hat\rho(k),ik_i\hat\rho(k)\rangle 
\nonumber\\ 
\nonumber\\ 
&=&\langle\fr{ik_j\hat\rho(k)}{k^2+m^2+(i|v|k_1+\lam)^2 
},ik_i\hat\rho(k)\rangle= 
\int\fr{k_ik_j|\hat\rho(k)|^2dk}{k^2+m^2+(i|v|k_1+\lam)^2}. 
\eeqn 
The matrix $H$ is well defined for $\Re\lam>0$ since the denominator 
does not vanish (or $g_\lam(x)$ exponentially decays). 
 The matrix $H$ is diagonal 
similarly to $K$. Indeed, 
if $i\ne j$, then at least one of these indexes is not equal to one, 
and the integrand in (\re{Cij}) 
is odd with respect to the corresponding variable.

As the result, $-\langle\na\Psi_1,\rho\rangle=HQ$, where $H$ 
 is the diagonal matrix with matrix elements $H_{jj}$, $1\le j\le3$. 
 Finally, the 
equations (\re{lte}) become 
\be\la{Mlam} 
M(\lam)\left( 
\ba{c} 
Q \\ P 
\ea 
\right)=\left( 
\ba{c} 
Q_0 \\ P_0' 
\ea 
\right),\,\,{\rm where}\,\,M(\lam)=\left( 
\ba{cc} 
\lam E  & -B_{v} \\ 
K-H(\lam) & \lam E 
\ea 
\right), 
\ee 
where the matrices  $K$ and $H(\lam)$ are diagonal. 
\\ 
{\it Step iv)} Assume for a moment that the matrix  
$M(\lam)$ is invertible for $\Re\lam>0$ 
 (later we shall prove that this the case indeed).  
Then  
\be\la{QP1} 
\left( 
\ba{c} 
Q \\ P 
\ea 
\right)=M^{-1}(\lam)\left( 
\ba{c} 
Q_0 \\ P_0' 
\ea 
\right),~~~~~~~~~\Re\lam>0. 
\ee 
Finally, formulas (\re{QP1})  
and (\re{Psi}) 
give the expression of the resolvent 
 $R(\lam)=(A-\lam)^{-1}$, $\Re\lam>0$.

\bl\la{cmf} 
The matrix function $M(\lam)$ ($M^{-1}(\lam)$) 
admits an analytic 
(meromorphic) continuation from the domain 
 $\Re\lam>0$ 
to the Riemann surface of the 
function $\sqrt{\lam^2+\mu^2}$. 
\el 
\pru 
The analytic continuation of  $M(\lam)$ exists by 
Lemma \re{cac} ii) and the convolution 
expressions in (\re{Cij}) since the function $\rho(x)$ 
is compactly supported by  (\re{ro}). 
 The inverse matrix is then meromorphic since it exists 
for large $\Re\lam$: this follows from (\re{Mlam}) 
since 
$H(\lam)\to 0$ as $\Re\lam\to\infty$ 
by  (\re{Cij}).\bo

\setcounter{equation}{0} 
\section{Analyticity in the Half-Plane} 
Here we prove the following proposition. 
\begin{pro}\la{analyt} 
The operator-valued function $R(\lam):{\cal E}\to{\cal E}$ 
is analytic for $\Re\lam>0$. 
\end{pro} 
{\bf Proof } 
It suffices to prove that the operator $A-\lam:{\cal E}\to 
{\cal E}$ has bounded inverse operator for $\Re\lam>0$. 
Recall that  $A=A_{v,v}$ where $|v|<1$. 
\medskip \\ 
{\it Step i)} Let us prove that Ker$\,(A-\lam )=0$ for $\Re\lam>0$. 
Indeed, assume that the vector 
$X_{\lam}=(\Psi_{\lam},\Pi_{\lam},Q_{\lam},P_{\lam})\in \cE$  
satisfies the equation 
 $(A-\lam)X_{\lam}=0$, 
that is $X_{\lam}$ is a solution to  
(\re{eq1}) with $\Psi_0=\Pi_0=0$ and 
$Q_0=P_0=0$. 
We must prove that $X_{\lam}=0$. 
 
Let us first show that $P_{\lam}=0$. 
Indeed, the trajectory 
$X:=X_{\lam}e^{\lam t}\in C(\R,\cE)$ 
is the solution to the equation $\dot X=AX$ of type  
(\re{line}) with $w=v$. 
Then 
$\cH_{v,v}(X(t))$ grows exponentially by  (\re{H0vv}), since the matrix $B_v$ is positive. 
This growth contradicts to the 
conservation of $\cH_{v,v}$, 
which 
follows from  Lemma \re{haml} ii) 
because $X(t)\in C^1(\R,\cE^+)$. 
The latter inclusion follows from 
Lemma \re{cres} 
since  $(\Psi_\lam,\Pi_\lam)$ satisfies equations (\re{Psi}) 
with $\Psi_0=\Pi_0=0$ and $Q=Q_\lam$. 
 
We now have $\lam Q_\lam=B_vP_\lam=0$ 
by the third equation of (\re{eq1}), and hence 
$Q_{\lam}=0$ because $\lam\ne 0$. 
Finally, 
 $\Psi_{\lam}=0$, $\Pi_{\lam}=0$ 
by equations (\re{Psi}) with $Q=Q_\lam=0$. 
\medskip \\ 
{\it Step ii)} One has 
$$ 
(A-\lam )\left( 
\ba{c} 
\Psi \\ \Pi \\ Q \\ P 
\ea 
\right)=\left( 
\ba{cccc} 
v\cdot\na-\lam & 1 & 0 & 0 \\ 
\De-m^2 & v\cdot\na-\lam & \cdot\na\rho & 0 \\ 
0 & 0 & -\lam  & B_{v} \\ 
\langle\cdot,\na\rho\rangle & 0 & \langle\na\psi_{v}, 
\cdot\na\rho\rangle & -\lam 
\ea 
\right)\left( 
\ba{c} 
\Psi \\ \Pi \\ Q \\ P 
\ea 
\right). 
$$ 
Thus, $A-\lam =A_0+T$, where 
$$ 
A_0=\left( 
\ba{cccc} 
v\cdot\na-\lam & 1 & 0 & 0 \\ 
\De-m^2 & v\cdot\na-\lam & 0 & 0 \\ 
0 & 0 & -\lam  & 0 \\ 
0 & 0 & 0 & -\lam 
\ea 
\right),\,\,\,~~~~~~T=\left( 
\ba{cccc} 
0 & 0 & 0 & 0 \\ 
0 & 0 & \cdot\na\rho & 0 \\ 
0 & 0 & 0 & B_{v} \\ 
\langle\cdot,\na\rho\rangle & 0 & \langle\na\psi_{v}, 
\cdot\na\rho\rangle & 0 
\ea 
\right). 
$$ 
The 
operator $T$ is finite-dimensional, and 
the operator $A_0^{-1}$ is bounded 
on $\cE$ by Lemma \re{cres}. 
Finally, $A-\lam =A_0(I+A_0^{-1}T)$, where 
$A_0^{-1}T$ is a compact operator. 
Since we know that Ker$\,(I+A_0^{-1}T)=0$, the operator $(I+A_0^{-1}T)$ 
is invertible 
by the Fredholm theory. 
 \bo 
 
\begin{cor} 
The matrix $M(\lam)$ of (\re{Mlam}) is invertible for $\Re\lam>0$. 
\end{cor} 
 
\setcounter{equation}{0} 
 
\section{Regularity on the Imaginary Axis} 
 
Next step should be an investigation of 
the limit values of the resolvent $R(\lam)$ 
at the imaginary axis $\lam=i\om$, $\om\in\R$, 
that is necessary for proving the decay (\re{Zdec}) 
of the solution $X(t)=(\Psi(t),\Pi(t),Q(t),P(t))$. 
 
Let us first 
describe  the continuous spectrum of the operator $A=A_{v,v}$ 
on the imaginary axis. 
By definition, 
the continuous spectrum corresponds to $\om\in\R$ 
such that the resolvent $R(i\om+0)$ is not a bounded  
operator on $\cE$. 
By the formulas (\re{Psi}), this is the case if 
 the Green function $g_\lam(y-y')$ fails to have 
the exponential decay. 
This is equivalent to the condition that 
$\Re\ka =0$, where $\ka$ is given by 
(\re{kappa}): $\ka=\ga\sqrt{\mu^2-\om^2}$. 
Thus, $i\om$ belongs to the continuous spectrum if (cf. (\re{kappa})) 
$$ 
|\om|\ge \mu=m\sqrt{1-v^2}. 
$$ 
By Lemma \re{cmf}, 
the limit matrix 
\be\la{M} 
M(i\om):=M(i\om+0)=\left( 
\ba{cc} 
i\om E & -B_v \\ 
K\!\!-\!\!H(i\om+0) & i\om E 
\ea 
\right), ~~~~~~~~~~\om\in\R, 
\ee 
exists, and its entries are continuous 
functions of $\om\in\R$, smooth for $|\om|< \mu$ 
and $|\om|>\mu$. 
Recall that the point $\lam=0$ belongs to the discrete spectrum of 
the operator $A$ by Lemma \re{ceig} i), and hence 
$M(i\om+0)$ is (probably) not invertible either at $\om=0$. 
\bp\la{regi} 
Let  (\re{ro}) and (\re{W}) hold, and 
$|v|<1$. Then 
the limit matrix 
 $M(i\om+0)$ is invertible for 
$\om\ne 0$, $\om\in\R$. 
\ep 
{\bf Proof } Let us consider   
the three possible cases 
$0<|\om|< \mu$, $|\om|= \mu$, 
and $|\om|>\mu$ separately. 
Let us remind that 
the matrices $K$ and $H$ are diagonal with the entries 
\be\la{al} 
K_{jj}=\int\fr{k_{j}^2|\hat\rho(k)|^2dk}{\ds k^2+m^2-(|v|k_1)^2}\,, 
\ee 
\be\la{Cii} 
H_{jj}(\lam)=\int\fr{k_j^2|\hat\rho(k)|^2dk}{k^2+m^2+(i|v|k_1+\lam)^2}, 
\,\,\,~~~~~~~~~~~~\Re\lam>0, 
\ee 
and $H_{22}=H_{33}$. 
Since 
 $v=(|v|,0,0)$,  the matrix $B_v$ is also diagonal: 
\be\la{Bv} 
B_v:=\nu(E-v\otimes v)= 
\left( 
\ba{ccc} 
\nu^3&0&0\\ 
0&\nu&0\\ 
0&0&\nu 
\ea 
\right) 
\ee 
since $\nu^2:=1-v^2$. 
Let us denote $F(\om):=-K+H(i\om+0)$ which is also 
diagonal, and let 
$F_{\Vert}:=F_{11}(\om)$, and $F_{\bot}:=F_{22}(\om)=F_{33}(\om)$. 
Then by (\re{M}) 
\beqn\la{detM} 
{\rm det}\,M(i\om)={\rm det}\,\left( 
\ba{cc} 
i\om E  & -B_v \\ 
-F(\om) & i\om E 
\ea 
\right)=-(\om^2+\nu^3F_{\Vert})(\om^2+\nu F_{\bot})^2,~~~\om\in\R. 
\eeqn 
The formula for the determinant 
is obvious since both matrices $F(\om)$ and $B_v$ are diagonal, hence the matrix 
$M(i\om)$ is equivalent to three independent  matrices $2\times2$. 
Namely, let us transpose the columns and rows of the matrix $M(i\om)$ 
in the order $(142536)$. Then we get the matrix with three $2\times 2$ 
blocks on the main diagonal. 
Therefore, the determinant of $M(i\om)$ is simply product of the determinants 
of the three matrices.

\noindent 
{\bf I.} First, let us consider the case 
$0<|\om|< \mu$. 
Then the invertibility  of $M(i\om)$ 
follows from  (\re{detM}) 
by the following lemma. 
\bl\la{lnW} 
For $0<|\om|<\mu$, the matrix $F(\om)$ 
is positive definite, i.e. $F_{jj}(\om)>0$, $j=1,2,3$. 
\el 
\pru 
First, let us check that the denominator in (\re{Cii}) 
is positive for $\lam=i\om$ with $|\om|<\mu$. 
 Indeed, it equals $m^2+k^2-(\om+|v|k_1)^2$ and we have to prove that 
$m^2+k^2>\om^2+2\om |v|k_1+v^2k_1^2$. By the condition 
 $|\om|<\mu=m\sqrt{1-v^2}$, it suffices to prove that 
$m^2+k^2\ge m^2(1-v^2)+2\om |v|k_1+v^2k_1^2$. This is equivalent 
to $k_2^2+k_3^2+m^2v^2+k_1^2(1-v^2)\ge2m|v|k_1\sqrt{1-v^2}$, 
which is evidently true. 
Thus, 
$$ 
F_{jj}(\om)=\int k_j^2|\hat\rho(k)|^2dk\left(\fr1{m^2+k^2-(|v|k_1+\om)^2}- 
\fr1{m^2+k^2-(|v|k_1)^2}\right),\,\,\,j=1,2,3. 
$$ 
 Let us prove that $F_{jj}(\om)>0$. Indeed, 
since $\hat\rho(k)=\hat\rho(-k)$, we obtain that 
\beqn\la{rep} 
\!\!\!\!\!\!\!\!\!\!\!\!\!\!\!\!F_{jj}(\om) 
=\int dk_2dk_3\int_0^{+\infty}k_j^2|\hat\rho(k)|^2\Bigg(\fr1{m^2+k^2- 
(|v|k_1+\om)^2}&+&\fr1{m^2+k^2-(|v|k_1-\om)^2} 
\nonumber\\ 
\nonumber\\ 
\!\!\!\!\!\!\!\!\!\!\!&-&\fr2{m^2+ 
k^2-(|v|k_1)^2}~\Bigg)dk_1. 
\eeqn 
Now it suffices to prove that the expression in brackets is positive 
(or positive infinite) 
under the conditions 
\be\la{cond} 
~~~~~~~|v|<1,~~~~~~~~0<|\om|\le \mu=m\sqrt{1-v^2}. 
\ee 
 This is proved 
in Appendix B. \bo 
\medskip 
 
\noindent 
{\bf II.} $\om=\pm\mu$. For example consider the case 
$\om=\mu$. 
Then formula  
(\re{Cii}) reads (see (\re{Dhat})): 
$$ 
H_{jj}(i\mu)=\int\fr{k_j^2|\hat\rho(k)|^2dk}{k_2^2+k_3^2+ 
(\nu k_1-m|v|)^2}. 
$$ 
Now the integrand   
has a unique singular point. The 
singularity is integrable, and therefore the terms 
$F_{jj}(\mu)$ are finite. Furthermore, the terms  
are positive  
by 
the integral representation 
(\re{rep}) again. Hence, 
 the matrix $M(i\mu)$ is invertible. 
\smallskip 
 
\noindent 
{\bf III.} $|\om|>\mu$. Here we apply another arguments: 
the invertibility  of $M(i\om)$ 
follows from  (\re{detM}) 
by  
the methods used in 
\ci[Chapter VII, formula (58)]{Vai}. 
\bl\la{lW} 
If (\re{W}) holds and if $|\om|>\mu$, 
then the imaginary part of the matrix $\ds\fr{\om}{|\om|}F(\om)$ is negative definite, i.e. 
 $\ds\fr{\om}{|\om|}\Im F_{jj}(\om)<0$, $j=1,2,3$. 
\el 
\pru 
Since $F(\om)=-K+H(i\om+0)$ where the matrix $K$ is real, 
it suffices to study the matrix $H(i\om+0)$. 
For $\ve>0$, we have 
\be\la{Hjjlim} 
H_{jj}(i\om+\ve)=\int\fr{k_j^2|\hat\rho(k)|^2dk}{k_1^2+k_2^2+k_3^2 
-(|v|k_1+\om-i\ve)^2+m^2},\,\,\,j=1,2,3. 
\ee 
Consider the denominator 
$$ 
\hat D_{\ve}(k)=k^2+m^2-(|v|k_1+\om-i\ve)^2. 
$$ 
It was shown above 
that $\hat D_{0}(k)\ne0$ if $|\om|<\mu$, and  
$\hat D_{0}(k)$ vanishes at one 
point if $|\om|=\mu$.  
On the other hand, for $|\om|>\mu$ 
the denominator $\hat D_{0}(k)$ vanishes 
 on the ellipsoid  
$$ 
T_{\om}=\{k:(\nu k_1-\fr{|v|\om}{\nu})^2+k_2^2+k_3^2=R^2:=\fr{\om^2-\mu^2} 
{\nu^2}\}, 
$$ 
where $\nu=\sqrt{1-v^2}$.  
We shall show below that 
the Plemelj formula for $C^1$-functions 
implies that 
\be\la{ImHjj} 
\Im 
H_{jj}(i\om+0)=-\fr{\om}{|\om|}\pi\int_{T_{\om}}\fr{k_j^2|\hat\rho(k)|^2} 
{|\na\hat D_{0}(k)|}dS, 
\ee 
where $dS$ is the element of the surface area.  
This immediately implies 
the statement of the Lemma since 
the integrand in (\re{ImHjj}) 
is positive by the Wiener condition (\re{W}).

Let us justify (\re{ImHjj}) for $\om>\mu>0$ 
(the case $\om<-\mu<0$ can be treated similarly). 
Let $\zeta\in C_0^{\infty}(\R^3)$ 
be a nonnegative 
cut off function equal one when $|\hat D_{0}(k)|<\de$ and vanishing 
when $|\hat D_{0}(k)|>2\de$. We fix a small $\de$ and split the integral 
(\re{Hjjlim}) in two parts: with the factor $\zeta$ in the integrand and 
with the factor $1-\zeta$. The limit of the second term as $\ve\to0$ 
is real. Hence, we have to 
calculate the imaginary part 
only for 
\be\la{Djjde} 
H_{jj}^{(\de)}(i\om+0)=\lim_{\ve\to0}\int\zeta(k)\fr{k_j^2|\hat\rho(k)|^2dk} 
{\hat D_{\ve}(k)}. 
\ee 
Denote $a(k)=\sqrt{k^2+m^2}$ and  $b(k)=|v|k_1+\om$. Then 
\be\la{1Dve} 
\fr1{\hat D_{\ve}(k)}=\fr1{a^2-(b-i\ve)^2}=\fr1{2a(a-b+i\ve)}+\fr1{2a(a+b-i\ve)}. 
\ee 
Note that $\hat D_{0}(k)\ne0$ if $b(k)=0$. Thus, $b(k)\ne0$ on $T_{\om}$, and therefore $b(k)\ne0$ on the support of $\zeta$ if $\de\ll1$. Since $b(k)>0$ when $v=0$ ($\nu=1$), we get that $b(k)>0$ on the support of $\zeta$ for all $v$ with $|v|<1$. 
 
We split the integral in (\re{Djjde}) in two terms accordingly to (\re{1Dve}). Then the second term is real for $\ve=0$. Now it remains to 
calculate the imaginary part 
of $h(i\om+0)$, where 
\be\la{Fiom} 
h(i\om+\ve):=\int\zeta(k)\fr{k_j^2|\hat\rho(k)|^2dk}{2a(a-b+i\ve)}. 
\ee 
One can rewrite (\re{Fiom}) as the iterated integral: 
 over the surfaces $T_{\om}^\al=\{k\in\R^3: 
a(k)-b(k)=\al,~~|\hat D_0(k)|<\de\}$ and over $\al$. Then we get 
$$ 
h(i\om+\ve)=\int\fr{u(\al)}{\al+i\ve}d\al,\,\,\,u(\al)=\int_{T_{\om}^\al}\zeta(k)\fr{k_j^2|\hat\rho(k)|^2dS}{2a|\na(a-b)|}, 
$$ 
and therefore 
$$ 
\Im h(i\om+0)=-\pi u(0)=-\pi\int_{T_{\om}}\zeta(k)\fr{k_j^2 
|\hat\rho(k)|^2dS}{2a|\na(a-b)|}. 
$$ 
This implies (\re{ImHjj}), since $\hat D_{0}(k)=a^2-b^2$ and 
$|\na\hat D_{0}(k)|=2a|\na(a-b)|$ on $T_{\om}$.\hfill\bo 
 
\noindent This completes the proofs of the Lemma 
15.3 and the Proposition 15.1. 
 
\bc\la{creg} 
Proposition \re{regi} implies that 
the matrix $M^{-1}(i\om)$ is smooth with respect to 
 $\om\in\R$ 
outside the three points $\om=0,\pm \mu$. 
\ec 
\br\la{reW} 
{\rm 
The proof of the  Lemma \re{lW} is the unique point 
in the paper where 
the Wiener condition is indispensable. 
In Lemma \re{lnW} we use 
only that the coupling function $\rho(x)$ 
is not identically zero. 
} 
\er 
 
\setcounter{equation}{0} 
\section{Singular Spectral Points}

Recall that the formula (\re{QP1}) expresses the Fourier-Laplace 
transforms $\ti Q(\lam),\ti P(\lam)$. Hence, 
 the vector components of the solution 
are given by the 
Fourier integral 
\be\la{QP1i} 
\left( 
\ba{c} 
Q(t) \\ P(t) 
\ea 
\right)=\ds\fr 1{2\pi}\int e^{i\om t}M^{-1}(i\om+0)\left( 
\ba{c} 
Q_0 \\ P_0' 
\ea 
\right)d\om 
\ee 
which converges in the sense of distributions. 
It remains to prove the continuity and decay of the  
vector components. 
The Corollary \re{creg} by itself 
is insufficient to prove the convergence and 
decay of the integral. 
Namely, we need an additional information 
about 
the regularity of the matrix $M^{-1}(i\om)$ at the 
singular 
points  $\om=0,\pm \mu$ 
and about some bounds at $|\om|\to\infty$. 
We shall study the points separately. 
\medskip\\ 
{\bf I.} 
Consider first  
the points $\pm \mu$. 
\begin{lemma} \la{Pui} 
The matrix  $M^{-1}(i\om)$ admits the following Puiseux expansion 
in a neighborhood of $\pm \mu$: there exists an $\ve_\pm>0$ s.t. 
\be\la{mom} 
M^{-1}(i\om)=\sum_{k=0}^{\infty}R_k^\pm(\om\mp\mu)^{k/2},\,\,\, 
|\om\mp\mu|<\ve_\pm,~~~~\om\in\R. 
\ee 
\end{lemma} 
{\bf Proof } 
It suffices to prove a similar expansion  for 
 $M(i\om)$. Then  (\re{mom})  
holds for $M^{-1}(i\om)$ as well, since the matrices 
$M(\pm i\mu)$ are invertible. The asymptotics for $M(i\om)$ holds by 
the convolution representation in (\re{Cij}): 
\be\la{Cjjj} 
H_{ij}(\lam)=\langle g_{\lam}*\pa_j\rho,\pa_i\rho\rangle 
\ee 
since $g_{\lam}$ admits the corresponding Puiseux 
expansions by formula (\re{glam}).\bo 
\medskip\\ 
{\bf II.} 
Second, we study the asymptotic behavior of 
$M^{-1}(\lam )$ at 
infinity. Let us recall that $M^{-1}(\lam )$ was originally defined for 
$\Re\lam >0,$ but it admits a meromorphic continuation to the 
Riemann 
surface of the function $\sqrt{\lam ^{2}+\mu ^{2}}$ 
(see Lemma \ref{cmf}). 
 
The following proposition 
is a very particular case of a general fundamental theorem 
about the bound for the truncated resolvent on the 
continuous spectrum. 
The bound plays a crucial role in the study of the long-time 
asymptotics of general linear hyperbolic PDEs, \ci{Vai}. 
\bp 
\la{162} 
We can find  
a matrix $R_{0}$ and a matrix-function 
$R_{1}(\om )$  such 
that 
$$ 
M^{-1}(i\om)=\fr {R_0}\om +R_1(\om ),~~~|\om|\ge\mu+1,~~~~~~~~~\om \in \R, 
$$ 
where 
\be\label{min} 
|\pa_\om^k R_1(\om )|\leq \fr{C_k}{|\om |^2}, 
~~~~~~~~~~~~|\om|\ge\mu+1,~~~~~~~~~\om \in \R 
\ee 
for every $k=0,1,2,...,$ 
\ep 
\pru 
By the structure (\ref{M}) of the matrix 
$M(i\om )$ it suffices to prove the 
following estimate for the 
elements of the matrix $H(i\om ):=H(i\om+0)$: 
\be \la{minC} 
|\pa _\om^k H_{jj}(i\om )|\leq \fr{C_k} 
{|\om|}, ~~~~~~~~~\om \in \R,~~~|\om|\ge\mu+1,~~~~j=1,2,3. 
\ee 
Let us rewrite (\ref{Cjjj}) as 
\be\la{Cjjja} 
H_{ij}(\lam)=\langle D^{-1}(\lam)\pa_j\rho,\pa_i\rho\rangle,~~~~~~~~\Re\lam>0, 
\ee 
where $D(\lam)$ is the operator  (\ref{fso}), 
and $D^{-1}(\lam)$ is a bounded operator on $L^2(\R^3)$. 
Let us denote by $B_R$ the ball $\{x\in\R^3:|x|<R\}$. 
Estimate (\ref{minC}) immediately follows 
from a more general bound 
\be \la{minCg} 
\Vert 
\pa _\om^k 
D^{-1}(i\om+0)f\Vert_{L^2(B_R)}\le \fr {C_k(R)}{|\om|} 
\Vert f\Vert_{L^2(B_R)},~~~~~~~~~\om \in \R,~~~|\om|\ge \mu+1 
\ee 
which holds for every $R>0$ and all functions 
$f(y)\in L^2_R:=\{f(y)\in L^2(\R^3): 
\supp f\subset B_R\}$. 
Namely,  
by (\re{ro}) the asymptotics 
(\ref{minC}) follows from the bound  
(\ref{minCg}) 
applied to the function $f(y)=\pa_j\rho(y)\in  L^2_R$ 
 with  $R\ge R_\rho$. 
The bound (\ref{minCg}) follows from a general 
estimate \ci[Thm 3]{Vai69} 
(see also \ci[the bound (A.2')]{Ag} 
\ci[Thm 8.1]{JK}, \ci[Thm 3]{Vai75}). 
\medskip\\ 
{\bf III.} 
Finally, consider the point  $\om=0$ which is 
the most singular. This is an isolated pole of a finite 
degree by Lemma \re{cmf}, and hence 
the Laurent expansion holds, 
\be\la{Lor} 
M^{-1}(i\om)=\sum_{k=0}^n L_k\om^{-k-1}+h(\om),~~~~~~~|\om|<\ve_0, 
\ee 
where $L_k$ are $6\times 6$ complex matrices, 
$\ve_0>0$, and $h(\om)$ is an analytic matrix-valued 
function for complex $\om$ with $|\om|<\ve_0$.

\section{Time Decay of the Vector Components} 
\setcounter{equation}{0} 
Here we prove the decay (\re{Zdec}) for the components $Q(t)$ and $P(t)$. 
\bl\la{171} 
Let $X_0\in \cZ_{v,\beta}$. Then $Q(t)$, $P(t)$ are continuous and 
the following bound holds, 
\be\la{decQP} 
|Q(t)|+|P(t)|\le \ds\fr {C(\rho, \ov v,d_0)}{(1+|t|)^{3/2}}, 
~~~~~~~t\ge 0. 
\ee 
\el 
\pru 
Expansions (\re{mom}), 
(\re{min}),  and (\re{Lor}) imply 
the convergence of the Fourier integral (\re{QP1i}) 
in the sense of distributions to a continuous function 
of $t\ge 0$. Let us prove  
the decay 
(\re{decQP}). 
We 
know that the linearized dynamics admits the secular  
solutions without decay, 
see (\re{secs}). 
The formulas (\re{inb}) give 
the corresponding components $Q_S(t)$ and $P_S(t)$ of the 
secular solutions, 
\be\la{secQPs} 
\left(\ba{c}Q_S(t)\\P_S(t)\ea\right) 
= 
\sum_1^3 C_j 
\left(\ba{c}e_j\\0\ea\right) 
+ 
\sum_1^3 D_j 
\Bigg[\left(\ba{c}e_j\\0\ea\right)t 
+\left(\ba{c}0\\ \pa_{v_j}p_v\ea\right)\Bigg]. 
\ee 
We claim that  the 
symplectic orthogonality condition leads to (\re{decQP}). 
Let us split the Fourier integral 
(\re{QP1i}) into three terms by  
using the partition of unity 
$\zeta_1(\om)+\zeta_2(\om)+\zeta_3(\om)=1$, $\om\in\R$: 
\beqn\la{QP1i3} 
\left( 
\ba{c} 
Q(t) \\ P(t) 
\ea 
\right)&=&\ds\fr 1{2\pi}\int e^{i\om t}(\zeta_1(\om)+\zeta_2(\om)+\zeta_3(\om))M^{-1}(i\om+0)\left( 
\ba{c} 
Q_0 \\ P_0' 
\ea 
\right)d\om 
\nonumber\\ 
\nonumber\\ 
&=&I_1(t)+I_2(t)+I_3(t), 
\eeqn 
where the functions $\zeta_k(\om)\in C^\infty(\R)$ are supported by 
\be\la{zsup} 
\left.\ba{rcl} 
\supp \zeta_1&\subset& \{\om\in\R:\ve_0/2\le|\om|\le\mu+2\} 
\\ 
\\ 
\supp \zeta_2&\subset& \{\om\in\R:|\om|\ge\mu+1\} 
\\ 
\\ 
\supp \zeta_3&\subset& \{\om\in\R:|\om|\le 
\ve_0\} 
\ea\right| 
\ee 
Then 
\smallskip\\ 
i) The function $I_1(t)\in C^\infty(\R)$ 
decays like $(1+|t|)^{-3/2}$ by the Puiseux 
expansion (\re{mom}). 
\\ 
ii) The function $I_2(t)\in C[0,\infty)$  
decays faster than any power of $t$ due to 
Proposition \re{162}. 
\\ 
iii) Finally, the function $I_3(t)$ generally does  
not decay if $n\ge 0$ 
in the Laurent expansion (\re{Lor}). 
Namely, the contribution of the analytic function $h(\om)$  
is a smooth function of $t\in\R$, and 
decays 
 faster than any power of $t$. On the other hand, the 
contribution of the Laurent series, 
\be\la{Fin} 
\left( 
\ba{c} 
Q_L(t) \\ P_L(t) 
\ea 
\right):= 
\ds\fr 1{2\pi} \int e^{i\om t}\zeta_3(\om) 
\sum_{k=0}^n L_k(\om-i0)^{-k-1} 
\left( 
\ba{c} 
Q_0 \\ P_0' 
\ea 
\right) 
d\om,~~~~~~~t\in\R, 
\ee 
is a polynomial function of $t\in\R$ 
of a degree $\le n$, 
modulo smooth  
functions of $t\in\R$ 
decaying faster than any power of $t$. 
This follows by the Cauchy theorem applied to the  
integral (\re{Fin}) if we change  
the integral over $\om\in [-\ve_0/2,\ve_0/2]$, 
where $\zeta_3(\om)\equiv 1$, 
 by the integral over the  
semicircle $e^{i\theta}\ve_0/2$, $\theta\in[\pi,0]$. 
Let us note that 
the formula (\re{secQPs}) 
gives an example of polynomial function arising from 
 (\re{Fin}).

We must show that 
the symplectic orthogonality condition  
eliminates the polynomial functions. 
Our main difficulty is that we know nothing about 
the order $n$ of the pole and about  
the Laurent coefficients 
$L_k$ 
of the matrix $M^{-1}(i\om)$ at $\om=0$. 
Our crucial observation has the following form: 
\smallskip\\ 
a) The components  (\re{secQPs}) 
of the secular solutions form a linear space $\cL_S$ 
of dimension dim~$\cL_S= 6$. 
\\ 
b) The polynomial functions  in 
 (\re{Fin}) belong to a linear space $\cL_L$ 
of dimension dim~$\cL_L\le 6$ since $(Q_0,P_0')\in\R^6$. 
\\ 
c) 
$\cL_S\subset\cL_L$ since any function (\re{secQPs}) 
admits a representation  
of the form 
(\re{Fin}). The validity of this representation 
 follows from 
the fact that the secular solutions  (\re{secs}) can be reproduced by 
our calculations with the Laplace transform. 
\smallskip\\ 
Therefore, we can conclude that 
\be\la{LLL} 
\cL_L =\cL_S. 
\ee 
Let us show that the secular 
solutions are forbidden since $X_0\in \cZ_{v,\beta}$, and 
hence the polynomial terms in (\re{Fin}) vanish, which implies 
the decay (\re{decQP}). 
 
First, the constructed vector components 
$Q(t)$ and $P(t)$ are continuous functions of $t\ge 0$.  
Hence,  
the corresponding field components  
$\Psi(t)$ and $\Pi(t)$ 
can be constructed 
by solving the first two equations of 
 (\re{Avv}), where $A_1$ is given by (\re{AA}) with $w=v=v(t_1)$ 
(see (\re{linf}) below). 
Proposition \re{pfi} i) in the next section implies that  
$X(t)\in C(\R,\cE)$. 
 
Second, 
the condition  $X_0\in \cZ_{v,\beta}$ implies that 
the entire trajectory $X(t)$ lies in  $\cZ_{v,\beta}$. This follows 
from the invariance of the space $\cZ_{v,\beta}$ under the generator 
$A_{v,v}$ (cf. Remark \re{rZ}). 
In other words, 
$X(t)=\bP_vX(t)$. 
 
On the other hand, 
 identity (\re{LLL}) implies that 
$X(t)$ can be corrected by a secular 
solution $X_S(t)$ s.t. the corresponding components  
$Q_\De(t)$ and $P_\De(t)$ 
of the difference $\De(t):=X(t)-X_S(t)$ decay at the rate  
 $(1+|t|)^{-3/2}$. 
Note that $\bP_v\De(t)=\bP_v X(t)=X(t)$ since 
$\bP_vX_S(t)=0$. 
 
Further, the difference $\De(t)\in C(\R,\cE)$ 
is a solution to  
the linearized equation (\re{Avv}). 
Hence, the corresponding   
norms of the  
field components of  
 $\De(t)$ also decay like  $(1+|t|)^{-3/2}$ that follows from  
Proposition \re{pfi} ii). 
Therefore, $\Vert \De(t)\Vert_{-\beta}\le C(1+|t|)^{-3/2}$, 
hence 
the components  $Q(t)$ and $P(t)$, 
of $X(t)=\bP_v \De(t)$ also decay like  $(1+|t|)^{-3/2}$.\bo

 
 
\setcounter{equation}{0} 
 
\section{Time Decay of Fields} 
In Sections 12-17 we denote by $X(t)$ the solution 
to the linearized equation  
 (\re{Avv}) with a fixed initial condition $X_0$. 
Here we  
consider an arbitrary solution 
$X(t)=(\Psi(\cdot,t),\Pi(\cdot,t),Q(t),P(t))$  
of the  
linearized equation. 
We shall prove a proposition which can be applied to  
the solution $X(t)$ from  previous 
sections as well as to the solution $\De(t)$ above.  
Let us study the field part of the solution, 
$ 
F(t)=(\Psi(\cdot,t),\Pi(\cdot,t)) 
$, 
solving the first two equations  from the system (\re{Avv}). 
These two equations have the form 
\be\la{linf} 
\dot F(t)=\left( 
\ba{ll} 
v\cdot\na & 1 \\ 
\De-m^2 & v\cdot\na 
\ea 
\right)F(t)+\left( 
\ba{l} 
0 \\ 
Q(t)\cdot\na\rho 
\ea 
\right). 
\ee 
We shall assume that the vector components decay, 
\be\la{linb} 
|Q(t)|\le \ds\fr {C(\rho, \ov v,d_0)}{(1+|t|)^{3/2}},~~~~ 
t\ge 0. 
\ee 
Proposition \re{lindecay} 
is reduced now to the following assertion. 
\begin{pro} \la{pfi} 
i) Let $Q(t)\in C([0,\infty);\R^3)$, 
and $F_0\in\cF$. 
Then 
equation (\re{linf}) admits a unique 
solution $F(t)\in C([0,\infty);\cF)$ with the initial condition 
$F(0)=F_0$. 
\\ 
ii) If $F_0\in\cF_\beta$ and if the decay (\re{linb}) holds, 
then the corresponding fields also decay 
uniformly with respect to $v$: 
\be\la{lins} 
\Vert F(t)\Vert_{-\beta}\le \ds\fr{C(\rho, \ov v, \ti v, 
d_0,\Vert F_0\Vert_{\beta})}{(1+|t|)^{3/2}} 
,~~~~t\ge 0, 
\ee 
for $|v|\le \ti v$ with any $\ti v\in (0,1)$. 
\end{pro} 
{\bf Proof } 
{\it Step i)} 
The 
statement i) follows from 
the Duhamel representation 
\be\la{Duh} 
F(t)=W(t)F_0+\left[\int_0^tW(t-s)\left( 
\ba{l} 
0 \\ Q(s)\cdot\na\rho 
\ea 
\right)ds 
\right],~~~~~~t\ge 0, 
\ee 
where $W(t)$ is the dynamical group 
of the modified Klein-Gordon equation 
\be\la{184st} 
\dot F(t)=\left( 
\ba{cc} 
v\cdot\na & 1 \\ 
\De-m^2 & v\cdot\na 
\ea 
\right)F(t). 
\ee 
The group $W(t)$ can be expressed through the group $W_0(t)$ 
of the standard Klein-Gordon 
 equation 
\be\la{KG} 
\dot \Phi(t)=\left( 
~~ 
\ba{cc} 
0 & 1 \\ 
\De-m^2 ~&~ 0 
\ea 
~~ 
\right)\Phi(t). 
\ee 
Namely, 
the problem (\re{KG}) 
corresponds to (\re{184st}), when $v=0$, and it is easy to see that 
\be\la{shifted} 
[W(t)F(0)](x)=[W_0(t)F(0)](x+vt),~~~~x\in\R^3,~~t\in\R. 
\ee 
Denote by $W(x-y,t)$ and $W_0(x-y,t)$ the (distribution) 
integral matrix kernels of 
 the operators $W(t)$ and $W_0(t)$ respectively. Then (\re{shifted}) implies 
that 
\be\la{shifted1} 
W(x-y,t)=W_0(x-y+vt,t),~~~x,y\in\R^3,~~t\in\R. 
\ee 
The identity 
(\re{shifted}) 
implies also 
the energy 
 conservation law for the group $W(t)$. 
Namely, for $(\Psi(\cdot,t),\Pi(\cdot,t))=W(t)F(0)$ we have 
$$ 
\int[|\Pi(x,t)-v\cdot\na\Psi(x,t)|^2+|\na\Psi(x,t)|^2+m^2|\Psi(x,t)|^2]dx=\co, 
~~~~~~t\in\R. 
$$ 
In particular, this gives  that 
\begin{equation} 
\Vert W(t)F_{0}\Vert _{\mathcal{F}}\leq C(\overline{v})\Vert F_{0}\Vert _{%
\mathcal{F}},\,\,\,t\in \mathrm{I\kern-.1567emR}.  \label{E} 
\end{equation} 
This estimate and (\ref{Duh}) imply the statement i). \newline 
\textit{Step ii)} The statement ii) follows from the Duhamel representation (%
\ref{Duh}) and the next lemma. 
 
\bl\label{lemma} For any $\beta >3/2$, $\overline{v}<1$ and 
$F_{0}\in \mathcal{F}_{\beta }$, the following decay holds, 
\begin{equation} 
\Vert W(t)F_{0}\Vert _{-\beta }\leq \frac{C(\beta ,\overline{v})}{(1+t)^{3/2}%
}\Vert F_{0}\Vert _{\beta },~~~~~~t\geq 0,  \label{dede} 
\end{equation} 
for the dynamical group $W(t)$ corresponding to the modified 
Klein-Gordon equation (\ref{184st}) with $|v|<\overline{v}$. 
\el 
\pru 
The lemma can be proved by general methods of Jensen and Kato \cite{JK} 
relying on the fundamental Agmon's estimate \cite[the bound 
(A.2')]{Ag}. We give an independent short proof for the 
convenience of the reader. 
\smallskip 
\\ 
\textit{Step i)}  
The matrix kernel $W_{0}(x-y,t)$ of the group 
$W_{0}(t)$ can be written explicitly since the solution to 
(\ref{KG}) has the form (see \cite{EKS}) 
\begin{equation} 
\Psi (\cdot ,t)= 
[\frac{\partial } 
{\partial t}R(t)\ast \Psi 
_{0}+R(t)\ast \Pi _{0}],\,\,\,\Pi (\cdot ,t)=\dot{\Psi}(\cdot ,t). 
\label{KGsol} 
\end{equation} 
Here $R(t)=R(\cdot ,t)=R_{0}(\cdot ,t)+R_{m}(\cdot ,t)$, and 
\[ 
R_{0}(x,t)=\frac{\delta (t-|x|)}{4\pi t},\,\,\, 
R_{m}(x,t)= 
-\frac{m}{4\pi }\frac{J_{1}^{+}(m\sqrt{t^{2}-|x|^{2}})}{\sqrt{t^{2}-|x|^{2}}}, 
\] 
where 
\[ 
J_{1}^{+}(m\sqrt{s}):=\left\{ 
\begin{array}{ll} 
J_{1}(m\sqrt{s}), & s\geq 0 \\ 
0 & s<0, 
\end{array} 
\right. 
\] 
and $J_{1}$ is the Bessel function of order 1. From here and well 
known asymptotics of the Bessel function it follows that 
\begin{eqnarray*} 
W_{0}(z,t) &=&0,~~~~~~|z|>t, \\ 
|\partial _{z}^{\alpha }W_{0}(z,t)| &\leq &C(\delta 
)(1+t)^{-3/2},~~~~|z|\leq (1-\delta )t, 
\end{eqnarray*} 
for $t\geq 1,$ $|\alpha |\leq 1$ and any $\delta >0.$ \ From the 
last two relations and (\ref{shifted1}) it follows that, for any 
$\overline{v}<1$ and $\varepsilon =\frac{1-\overline{v}}{2},$ the 
following estimates hold for 
the matrix kernel $W(z,t)$ of the group $W(t):$%
\begin{eqnarray} 
W(z,t) &=&0,~~~~~~|z|>(1+\overline{v})t,  \label{1a1} \\ 
|\partial _{z}^{\alpha }W(z,t)| &\leq 
&C(\overline{v})(1+t)^{-3/2},~~~~~~~|z|\leq \varepsilon 
t,~~~~~|\alpha |<1.  \label{1z1} 
\end{eqnarray} 
\textit{Step ii)} Let us fix an arbitrary $t\geq 1$, and split the 
initial function $F_{0}$ in two terms, $F_{0}=F_{0,t}^{\prime 
}+F_{0,t}^{\prime \prime }$ such that 
\begin{equation} 
\begin{array}{l} 
\Vert F_{0,t}^{\prime }\Vert _{\beta }+\Vert F_{0,t}^{\prime 
\prime }\Vert _{\beta }\leq C\Vert F_{0}\Vert _{\beta 
},~~~~~~~t\geq 1, 
\end{array} 
\label{FFn} 
\end{equation} 
and 
\begin{eqnarray} 
F_{0,t}^{\prime }(x)=0, &~~~~~&|x|>\frac{\varepsilon 
t}{2},\bigskip 
\label{F'} \\ 
F_{0,t}^{\prime \prime }(x)=0, &~~~~~&|x|<\frac{\varepsilon t}{4}, 
\label{F''} 
\end{eqnarray} 
where $\varepsilon >0$ is defined in (\ref{1z1}). The estimate for 
$W(t)F_{0,t}^{\prime \prime }$ follows by (\ref{E}), (\ref{F''}) 
and (\ref {FFn}): 
\begin{eqnarray} 
\Vert W(t)F_{0,t}^{\prime \prime }\Vert _{-\beta } &\leq &\Vert 
W(t)F_{0,t}^{\prime \prime }\Vert _{\mathcal{F}}\leq C\Vert 
F_{0,t}^{\prime 
\prime }\Vert _{\mathcal{F}}  \nonumber \\ 
&\leq &C_{1}(\varepsilon )\Vert F_{0,t}^{\prime \prime }\Vert 
_{\beta }(1+t)^{-\beta }\leq C_{2}(\varepsilon )\Vert F_{0}\Vert 
_{\beta }(1+t)^{-\beta },~~~~~~~~~~t\geq 1.  \label{zxc} 
\end{eqnarray} 
\textit{Step iii)} It remains to estimate $W(t)F_{0,t}^{\prime }$. 
We split the operator $W(t)$, for $t>1$, in two terms: 
\[ 
W(t)=(1-\zeta )W(t)+\zeta W(t), 
\] 
where $\zeta $ is the operator of multiplication by the function $\zeta (%
{|x|}/{t})$ such that $\zeta =\zeta (s)\in C_{0}^{\infty }(\mathrm{I\kern%
-.1567emR}),$ $\zeta (s)=1$ for $|s|<\varepsilon /4,$ $\zeta (s)=0$ for $%
|s|>\varepsilon /2.$ Since 
\[ 
|\partial _{x}^{\alpha }\zeta ({|x|}/{t})|\leq C,~~~~~ 
|\alpha |\leq 1,~~~~~~t\geq 1, 
\] 
and $1-\zeta ({|x|}/{t})=0$ for $|x|<\varepsilon t/4,$ we have, for $%
t\geq 1,$%
\[ 
||(1-\zeta )W(t)F_{0,t}^{\prime }||_{-\beta }\leq 
C_{3}(\varepsilon )(1+t)^{-\beta }||(1-\zeta )W(t)F_{0,t}^{\prime 
}||_{\mathcal{F}}\leq C_{4}(\varepsilon )(1+t)^{-\beta 
}||W(t)F_{0,t}^{\prime }||_{\mathcal{F}}. 
\] 
From here, (\ref{E}) and (\ref{FFn}) it follows that 
\begin{equation} 
||(1-\zeta )W(t)F_{0,t}^{\prime }||_{-\beta }\leq 
C_{5}(\varepsilon )(1+t)^{-\beta }||F_{0,t}^{\prime 
}||_{\mathcal{F}}\leq C_{6}(\varepsilon )(1+t)^{-\beta 
}||F_{0}||_{\mathcal{F}},~~~~~t\ge 1.  \label{zaz} 
\end{equation} 
\textit{Step iv)} 
Thus, in order to complete the proof of Lemma \ref{lemma}, it 
remains to receive a similar estimate for $\zeta 
W(t)F_{0,t}^{\prime }$. 
Let $\chi _{\varepsilon t/2}$ be the characteristic function of the ball $%
|x|\leq \varepsilon t/2.$ We will use the same notation for the 
operator of multiplication by this characteristic function. From 
(\ref{F'}) it follows that 
\[ 
\zeta W(t)F_{0,t}^{\prime }=\zeta W(t)\chi _{\varepsilon 
t/2}F_{0,t}^{\prime }. 
\] 
The matrix kernel $W^{\prime }(x,y,t)$ of the operator $\zeta 
W(t)\chi _{\varepsilon t/2}$ is equal to 
\[ 
W^{\prime }(x,y,t)=\zeta ({|x|}/{t})W(x-y,t)\chi _{\varepsilon 
t/2}(y). 
\] 
Since $\zeta ({|x|}/{t})=0$ for $|x|>\varepsilon t/2\ $and 
$\chi 
_{\varepsilon t/2}(y)=0$ for $|y|>\varepsilon t/2,$ the estimate 
(\ref{1z1})\ implies that 
\begin{equation} 
|\partial _{x}^{\alpha }W^{\prime }(x,y,t)|\leq C(\overline{v})(1+t)^{-3/2}, 
~~~~~|\alpha |<1,~~~~~~t\geq 1.  \label{qaz} 
\end{equation} 
The norm of the operator $\zeta W(t)\chi _{\varepsilon t/2}: 
\mathcal{F}_{\beta }\rightarrow \mathcal{F}_{-\beta }$ is equivalent to the 
norm of the operator 
\[ 
(1+|x|)^{-\beta }\zeta W(t)\chi _{\varepsilon t/2}(1+|y|)^{-\beta }: 
\mathcal{F}\rightarrow \mathcal{F}. 
\] 
The norm of the later operator does not exceed the sum in $\alpha 
,~|\alpha |\leq 1$ of the norms of operators 
\begin{equation} 
\partial _{x}^{\alpha }[(1+|x|)^{-\beta 
}\zeta W(t)\chi 
_{\varepsilon t/2}(1+|y|)^{-\beta }]: 
L^{2}(\mathrm{I\kern-.1567emR}^{3})\oplus L^{2}(\mathrm{I\kern-.1567emR}^{3}) 
\rightarrow L^{2}(\mathrm{I\kern-.1567emR}^{3}) 
\oplus L^{2}(\mathrm{I\kern-.1567emR}^{3}). 
\label{1234} 
\end{equation} 
From (\ref{qaz}) it follows that operators (\ref{1234}) are 
Hilbert-Schmidt operators since $\beta >3/2,$ and their 
Hilbert-Schmidt norms do not exceed $C(1+t)^{-3/2}$. Hence 
\begin{equation} 
||\zeta W(t)F_{0,t}^{\prime }||_{-\beta }\leq 
C(\overline{v})(1+t)^{-3/2}||F_{0,t}^{\prime }||_{\beta } 
\leq C_{7}(\overline{v})(1+t)^{-3/2}||F_{0,t}||_{\beta }, 
\,~~~\,\,t\geq 1. 
\label{HS} 
\end{equation} 
The last estimate above is due to (\ref{FFn}). Finally, the 
estimates (\ref {HS}), (\ref{zaz}) and (\ref{zxc}) imply 
(\ref{dede}). 
${\hfill \hbox{\enspace{\vrule height 7pt depth 0pt 
width 7pt}}}$

 

\appendix 
 
\setcounter{equation}{0} 
 
\section{Appendix: Computing Symplectic Form} 
 
Here we justify the formulas  (\re{Omega})-(\re{alpha}) 
for 
the matrix $\Om$. 
 For $j,l=1,2,3$ it follows from (\re{inb}) and (\re{OmJ}) 
that 
 \be\la{jl} 
\Om(\tau_j,\tau_l) 
=\langle\pa_j\psi_v,\pa_l\pi_v\rangle-\langle\pa_j\pi_v,\pa_l\psi_v\rangle, 
\ee 
\be\la{jp3lp3} 
\Om(\tau_{j+3},\tau_{l+3})=\langle\pa_{v_j}\psi_v,\pa_{v_l}\pi_v\rangle- 
\langle\pa_{v_j}\pi_v,\pa_{v_l}\psi_v\rangle, 
\ee 
and 
\be\la{Wx} 
\Om(\tau_{j},\tau_{l+3}) 
=-\langle\pa_{j}\psi_v,\pa_{v_l}\pi_v\rangle+ 
\langle\pa_{j}\pi_v,\pa_{v_l}\psi_v\rangle+e_j\cdot\pa_{v_l}p_v. 
\ee 
Let us transfer to the Fourier representation. Set 
\be\la{fur} 
\hat\psi(k):=(2\pi)^{-3/2}\int e^{ikx}\psi(x)dx. 
\ee 
It is easy to compute that 
\be\la{hpsiv} 
\hat\psi_v(k)=-\fr{\hat\rho(k)}{k^2+m^2-(kv)^2},\,\,\, 
\hat\pi_v(k)=i(kv)\hat\psi_v(k). 
\ee 
Further, differentiating, we obtain 
\be\la{derv} 
\pa_{v_j}\hat\psi_v=\fr{2(kv)k_{j}}{k^2+m^2-(kv)^2}\hat\psi_v,\,\,\, 
\pa_{v_j}\hat\pi_v=ik_{j}\fr{k^2+m^2+(kv)^2}{k^2+m^2- 
(kv)^2}\hat\psi_v,\,\,\,j=1,2,3, 
\ee 
and 
$$ 
\pa_{v_j}p_v:= 
\fr{e_j}{\sqrt{1-v^2}}+\fr{v_{j}}{(1-v^2)^{3/2}}v,\,\,\,j=1,2,3. 
$$ 
Then for $j,l=1,2,3$ we see from (\re{jl}) by the  Parseval identity that 
\be\la{jl123} 
\Om(\tau_j,\tau_l)=-2i\int\, k_jk_l(kv)|\hat\psi_v|^2\,dk 
= 0, 
\ee 
since the integrand is odd in $k$. Similarly, by (\re{jp3lp3}), 
\be\la{jl456} 
\Om(\tau_{j+3},\tau_{l+3})=-4i\int\fr{k_jk_l(kv)(k^2+m^2+(kv)^2)|\hat\psi_v|^2}{(k^2+m^2-(kv)^2)^2}=0. 
\ee 
Finally, by (\re{Wx}), 
\beqn\la{Wk} 
\Om(\tau_j,\tau_{l+3})&=&\int dk|\hat\psi_v|^2k_jk_{l}\left[\fr{k^2+ 
m^2+(kv)^2}{k^2+m^2-(kv)^2}+\fr{2(kv)^2} 
{k^2+m^2-(kv)^2}\right]+e_j\cdot\pa_{v_l}p_v 
\nonumber\\\nonumber\\ 
 &=&\int dk|\hat\psi_v|^2k_jk_{l}\fr{k^2+m^2+3(kv)^2}{k^2+m^2-(kv)^2}+ 
e_j\cdot\left(\fr{e_l}{\sqrt{1-v^2}} 
+\fr{v_{l}v}{(1-v^2)^{3/2}}\right). 
\eeqn 
This completes the proof of (\re{Omega}) - (\re{alpha}).

 
\section{Appendix: Positivity of the Matrix $F$} 
\setcounter{equation}{0} 

Here we justify the inequality used above 
in the proof of Lemma \re{lnW}: 
$$ 
\fr1{m^2+k^2-(|v|k_1+\om)^2}+\fr1{m^2+k^2-(|v|k_1-\om)^2}-\fr2{m^2+k^2- 
(|v|k_1)^2}>0 
$$ 
under the conditions (\re{cond}): 
\be\la{condB} 
~~~~~~~~|v|<1,\,\,\,\,\,\,0<|\om|\le \mu=m\sqrt{1-v^2}. 
\ee 
Let us denote  $M^2:=m^2+k^2$, 
$r_\pm:=|v|k_1\pm\om$, and $r:=|v|k_1$. Then 
the inequality reads, after cancellation by $2M$, 
\be\la{B2} 
\fr1{M-r_+}+\fr1{M-r_-}-\fr2{M-r} 
+ 
\fr1{M+r_+}+\fr1{M+r_-}-\fr2{M+r}>0. 
\ee 
The sum of the first three terms 
in (\re{B2}) 
can be written as 
\be\la{B3} 
\fr1{N-\om}+\fr1{N+\om}-\fr2{N}=  \fr{2\om^2}{(N+\om)(N-\om)N} 
\ee 
where $N:=M-r$. 
It is easy to check that $N\pm\om\ge 0$ 
and $N>0$ 
under the conditions (\re{condB}). Hence, the sum 
(\re{B3}) is positive (or positive infinite). 
Similarly, 
the sum of the last three terms in (\re{B2}) also is positive.



\end{document}